\documentclass[10pt,a4paper]{article}

\usepackage{geometry}
\usepackage[utf8]{inputenc}
\usepackage{amsthm, amssymb, natbib, graphicx}
\usepackage{listings}
\usepackage[ruled,vlined]{algorithm2e}
\usepackage{caption}
\usepackage{amsmath,tikz}

\usepackage{algpseudocode}
\usepackage{array, makecell}
\usepackage{sidecap, caption}
\usepackage{subcaption}
\usepackage{xcolor, color} 
\definecolor{darkred}{RGB}{139,0,0}
\definecolor{darkblue}{RGB}{0,0,139}
\usepackage{enumitem}
\usepackage{comment} 
\usepackage{algorithm2e}
\usepackage[most]{tcolorbox}
\setcitestyle{numbers, comma,square}
\usepackage{bbm}
\usepackage{caption}
\usepackage{tabularray}
\usepackage{diagbox}
\UseTblrLibrary{booktabs}
\usepackage{csquotes}
\usepackage[
    pdfpagelabels,
    hypertexnames=false,  
    bookmarksopen=true,
    bookmarksopenlevel=1,  
    colorlinks=true,
    pdfstartview=FitV,
    linkcolor=blue,
    citecolor=blue,
    urlcolor=blue,
]{hyperref} 
\hypersetup{
    colorlinks=true,
    linkcolor=blue,
    filecolor=magenta,      
    urlcolor=cyan,
    pdftitle={When is Polynomial Optimization SoS-Exact? Sufficient Regularization and Strict Gram Representations},
    pdfauthor={Coralia Cartis and Wenqi Zhu},
    pdfpagemode=FullScreen,
    }
\usepackage{cleveref}

\newtheorem{proposition}{Proposition}[section]
\newtheorem{theorem}{Theorem}[section]
\newtheorem{corollary}{Corollary}[section]
\newtheorem{lemma}{Lemma}[section]
\newtheorem{example}{Example}[section]
\newtheorem{remark}{Remark}[section]

\newtheorem{definition}{Definition}[section]

\usepackage[normalem]{ulem}

\newcommand{\ignore}[1]{}

\DeclareMathOperator*{\argmin}{argmin}

\DeclareMathOperator{\trace}{trace}

\newcommand{\R}{\mathbb{R}}

\newcommand{\SoS}{\mathrm{SoS}}
\newcommand{\coeff}{\operatorname{coeff}}

\newcommand{\lmin}{\lambda_{\min}}

\title{When is Polynomial Optimization SoS-Exact?\\ Sufficient Regularization and Strict Gram Representations}
\usepackage{authblk}
\author[1,$\dagger$]{Coralia Cartis}
\author[1,$*$]{Wenqi Zhu}
\affil[1]{%
Mathematical Institute, University of Oxford\par
Woodstock Road, Oxford OX2 6GG, United Kingdom
}
\affil[$\dagger$]{\texttt{coralia.cartis@maths.ox.ac.uk}}
\affil[$*$]{\texttt{wenqi.zhu@maths.ox.ac.uk}}
\date{}

\begin{document}
\maketitle

\begingroup
\renewcommand\thefootnote{}
\footnotetext{$*$ The authors are listed in alphabetical order. Wenqi Zhu is the
primary contributor.}
\addtocounter{footnote}{-1}
\endgroup

\begin{abstract} 
Global minimization of multivariate polynomials is computationally
difficult in general, and nonnegativity does not necessarily imply
a sum-of-squares (SoS) representation. We show that adding sufficiently
strong Euclidean regularization to a polynomial provides a general mechanism for
overcoming such challenges. Namely, for any
multivariate polynomial $P$ of degree at most $p$ with a nonzero
linear coefficient, and any even integer $2d>p$, we prove that
\[
m_\sigma(s)=P(s)+\frac{\sigma}{2d}\|s\|_2^{2d}
\]
is SoS exact for all sufficiently large $\sigma$.
Furthermore, the shifted polynomial
$m_\sigma-\min m_\sigma$ admits a positive-definite Gram
representation on appropriately-shifted nonconstant monomials. Consequently, the global minimum
and minimizer can be recovered exactly by a single degree-$d$
semidefinite programming (SDP) representation, whose size is
polynomial in the dimension for fixed $d$, without requiring
higher levels of the SoS hierarchy.
The proof develops a degree-independent mechanism based on
the Bregman remainder of the regularizer: its strict Gram
margin absorbs the lower-degree polynomial perturbations.
We further establish an analogous exact SDP regime for
polynomial optimization over sufficiently small Euclidean
balls, connecting large regularization with small-radius
trust-region problems. Finally, we investigate how the
algebraic structure of the regularizer governs the existence
of a strict Gram margin, extending the mechanism beyond
Euclidean norms while identifying limitations of
coordinate-separable regularization.
\end{abstract}

\section{\texorpdfstring{{Introduction}}{Introduction}}
\label{intro}

Global minimization of a multivariate polynomial is a fundamental
yet computationally challenging problem, even at relatively low degrees.
For instance, minimizing a homogeneous cubic over the Euclidean sphere
is NP-hard \cite{nesterov2003random}; testing nonnegativity of quartic
polynomials is NP-hard \cite{murty1987some}; and even deciding whether
a quartic polynomial has a local minimizer is strongly NP-hard
\cite{ahmadi2022complexity}.
Despite these general hardness results, certain polynomial classes
admit exact and computationally tractable formulations due to their
underlying algebraic structure.
This raises a fundamental question: beyond the classical tractable
cases, can one identify broader classes of multivariate polynomials
whose global minima can be computed and certified exactly?

A classical illustration of such algebraic structure dates back to
Hilbert's celebrated work in 1888. He showed that nonnegative
polynomials need not be sums of squares (SoS), and characterized
the three classical cases in which the two notions coincide:
univariate polynomials of arbitrary degree, quadratic polynomials
in any number of variables, and bivariate quartics
\cite{hilbert1888darstellung}.
This distinction has fundamental consequences for global optimization.
Indeed, if a polynomial $h$ attains its global minimum $h^*$,
then the shifted polynomial $h-h^*$ is nonnegative.
When this nonnegative polynomial admits an SoS representation,
its global minimum can be computed and certified through
semidefinite programming (SDP)
\cite{parrilo2000structured,lasserre2001global}.
Beyond Hilbert's classical cases, it is therefore natural to ask
whether additional algebraic structure can guarantee SoS exactness
and, consequently, exact global minimization by a tractable SDP.

In this paper, we show that sufficiently strong regularization
provides precisely such a mechanism. We consider a broad class
of multivariate polynomials of arbitrary even degree, obtained
by adding a higher-degree Euclidean regularizer to a polynomial
of lower degree. Let $P$ be an arbitrary multivariate polynomial
of degree at most $p$, and let $2d>p$ be an even integer.
We study the regularized polynomial
\begin{equation}
\label{eq:intro-regularized-model}
m_\sigma(s)
:=
P(s)+\frac{\sigma}{2d}\|s\|^{2d},
\qquad \sigma>0,
\end{equation}
where $\|\cdot\|$ denotes the Euclidean norm.
Since the regularization term has higher degree than $P$,
$m_\sigma$ is coercive and therefore always attains
a global minimum. 
Our main result shows that, under a nonzero linear-term
assumption, sufficiently strong regularization makes
$m_\sigma$ \emph{SoS exact}, with its unique global minimizer
recoverable from a single SDP of polynomial size
for every fixed $d$.
To be precise, let
\[
s_\sigma^*\in\arg\min_{u\in\R^n}m_\sigma(u),\qquad
m_\sigma^*:=m_\sigma(s_\sigma^*),\qquad
q_\sigma:=m_\sigma-m_\sigma^*,
\]
where $s_\sigma^*$ denotes a minimizer and  $m_\sigma^*$, its optimal value.
The shifted polynomial $q_\sigma$ is nonnegative by construction.
We say that $m_\sigma$ is \emph{SoS exact} if $q_\sigma$
is a sum of squares of polynomials of degree at most $d$.
Equivalently, the global minimum $m_\sigma^*$ is attained
exactly by a single degree-$d$ SoS representation, with no
higher-order relaxations required.
Furthermore, the SDP recovers not only the exact global minimum
but also its unique global minimizer. In particular, the
optimal moment matrix has rank one, and the minimizer can be
read directly from the moments corresponding to the linear
monomials. 
Our main result can be stated as follows; the precise
statements are given in
\Cref{cor:SoS-taylor-model,cor:polynomial-size-sdp}.

\begin{theorem}[Main result, informal]
\label{thm:intro-main}
Let $P$ be a polynomial of degree at most $p$ with
nonzero linear coefficient $g:=\nabla P(0)\neq0$,
let $2d>p$ be even, and let $m_\sigma$ be defined
in \eqref{eq:intro-regularized-model}.
Then there exists a threshold $\sigma^*>0$
such that, for every $\sigma\ge\sigma^*$:
\begin{enumerate}[label=\textup{(\roman*)}]
\item The shifted polynomial $m_\sigma-m_\sigma^*$
is a sum of squares of polynomials of degree at most $d$.
Furthermore, it admits a positive-definite Gram
representation on the nonconstant monomials
centred at global minimizer $s_\sigma^*$.

\item 
A single degree-$d$ moment SDP
attains the exact global minimum and has a unique
rank-one optimal moment matrix, from which $s_\sigma^*$
can be recovered directly. For fixed $d$, this SDP
has size polynomial in $n$.
\end{enumerate}
\end{theorem}

\noindent The threshold $\sigma^*$ is given explicitly
in terms of the coefficients of $P$.
Numerical evidence in the quartic case suggests
that exactness may occur at substantially smaller
regularization parameters
\cite{zhu2026sufficiently}.



The class in \eqref{eq:intro-regularized-model} arises naturally
in higher-order nonlinear optimization, where Taylor models
combined with regularization provide a standard globalization
mechanism. Any polynomial of degree at most $p$ can be written
exactly in Taylor form,
\begin{equation}
\label{eq:intro-taylor-model}
P(s)=\mathcal T_p(s)
:=\sum_{j=0}^p\frac1{j!}\,T_j[s]^j,
\end{equation}
where $T_j$ is a symmetric $j$th-order tensor determined by
the coefficients of $P$ (see \Cref{par:coefficient-matching}).
For a Taylor approximation of a smooth function $f$ at $x$,
we have $T_j=\nabla^j f(x)$. The associated regularized model is
\begin{equation}
\label{eq:intro-ARp-model}
m_\sigma(s)=\mathcal T_p(s)
+\frac{\sigma}{p'}\|s\|^{p'},
\end{equation}
where $\sigma>0$ and our polynomial setting takes
the even exponent $p'=2d>p$. Thus the $P$ in
\eqref{eq:intro-regularized-model} and the $\mathcal T_p$ here
are identical. 
Beginning with cubic regularization of Newton's method
\cite{nesterovpolyak2006}, higher-order adaptive regularization
methods have developed a rich convergence and complexity theory
\cite{birgin2017worst,cartis2019universal,cartis2022evaluation,
doikov2022local}. Yet the global solution of their polynomial
subproblems remains a central computational challenge.
Nesterov \cite{nesterov2021tensor} explicitly asked whether
tools from algebraic geometry, particularly sums of squares,
could be used to minimize such subproblems. Our result
provides a positive answer in the sufficiently regularized
regime. We also study the complementary trust-region
globalization strategy \cite{conn2000trust}, which minimizes
the same Taylor polynomial over $\|s\|\le\Delta$.
For sufficiently small radii, we establish a correspondence
with the large-regularization regime; see
\Cref{sec:TRplus-formulation}.

Classically, convexity is one of the principal
sources of tractability in optimization
\cite{rockafellar1993lagrange}.
We emphasize, however, that our framework does
not rely on convexifying the polynomial model.
Instead, we retain the regularized Taylor
polynomial itself and prove that sufficiently
strong Euclidean regularization makes its
shifted nonnegative polynomial SoS exact,
even when the model remains nonconvex.
Thus regularization can do more than provide
boundedness or globalization: it can enforce
algebraic exactness without requiring convexity. 
This provides a different route from the
SoS-convex higher-order Newton framework of
Ahmadi, Chaudhry and Zhang
\cite{ahmadi2024higherorder}.
Their approach constructs an SoS-convex
polynomial subproblem through suitable
corrections and regularization, whereas our
result establishes an exact SoS representation
for the original regularized model, without
assuming convexity or SoS-convexity.

In addition to the result in Theorem \ref{thm:intro-main}, our contributions reveal three complementary aspects of this algebraic exactness mechanism.

\paragraph{1. A mechanism linking optimality and algebraic geometry.}
Our central technical result shows that the Bregman remainder of the Euclidean
regularizer admits a positive-definite Gram representation on the full
vector of nonconstant monomials. When the regularization is sufficiently
large, this strict Gram margin absorbs the remaining Taylor terms,
yielding an exact SoS representation of the shifted nonnegative polynomial. 
This mechanism is distinct from our earlier quartic analysis
\cite{zhu2026sufficiently}, which relies on identities specific to degree
four. The Bregman/strict-Gram argument developed here applies at every
even regularization degree and extends to broader classes of regularizers;
see \Cref{sec: strict gram beyond}. Thus the present result is not merely
a higher-degree extension of the quartic case, but a degree-independent
approach to establishing SoS exactness.


\paragraph{2. Small-radius trust-region problems admit exact SDP solutions.}
The same mechanism governs polynomial optimization over a Euclidean ball.
Minimizing a degree-$p$ polynomial over $\|s\|\le\Delta$ is a
higher-order trust-region subproblem and is NP-hard in general
for $p\ge3$ \cite{nesterov2003random}. In
\Cref{sec:TRplus-formulation}, we show that when the linear coefficient
is nonzero and the radius $\Delta$ is sufficiently small, the global
minimizer is unique, lies on the sphere $\|s\|=\Delta$, and can be
recovered exactly from a rank-one solution of a degree-$d$ moment
relaxation (\Cref{thm:small-radius-TRplus}). 
At an active trust-region solution, the multiplier
$\lambda_\Delta$ determines the matched regularization parameter
$\sigma=\lambda_\Delta/\Delta^{2d-2}$.
As $\Delta\downarrow0$, this parameter tends to infinity,
connecting small trust regions with the large-regularization
regime of our main theorem.
The radius threshold is explicit in the tensor norms of the
Taylor coefficients. Thus, the two classical globalization
strategies share the same algebraic mechanism for SoS exactness.

\paragraph{3. The choice of regularization matters.}
The mechanism depends on the algebraic structure of the regularizer,
not merely on its strength. Euclidean regularization produces the
strict Gram margin underlying our SoS exactness result, whereas
other regularizers need not possess this property, regardless of
how large the regularization parameter becomes. In
\Cref{sec: strict gram beyond}, we investigate which regularizers
admit such a margin: weighted Euclidean norms and regularizers
constructed from spanning families of linear forms inherit the
strict Gram property, while coordinate-separable regularization,
although coercive and SoS-convex, does not control all mixed
monomial directions needed for a full positive-definite Gram
representation.
This distinction holds across the even regularization degrees
considered here. It shows that the choice of regularizer,
often treated as a modelling tool, plays a decisive
role in the algebraic structure underlying SoS exactness.


 \paragraph{Background and Notation}

We recall only the basic notions used throughout the paper: SoS
polynomials, monomial vectors, Gram representations, and SoS-convexity.
Unless stated otherwise, \(\|\cdot\|\) denotes the Euclidean norm.

\begin{definition}[Sum of squares and Gram representation]
\label{def:sos-gram}
Let \(h\in\R[x]\) be a polynomial of degree at most \(2d\).  We say that
\(h\) is a \emph{sum of squares} (SoS) if there exist polynomials
\(h_1,\ldots,h_r\) of degree at most \(d\) such that
$$
    h(x)=\sum_{i=1}^r h_i(x)^2.
$$
We denote by \(\Sigma[x]_{2d}\) the cone of SoS polynomials of degree at most
\(2d\). 
For \(k\ge0\), let \([x]_k\) denote the normalized vector of monomials of
total degree \(k\), chosen so that \(\|[x]_k\|^2=\|x\|^{2k}\).  We write
$
    z_d(x):=\bigl(1,[x]_1^\top,\ldots,[x]_d^\top\bigr)^\top
$
for the full monomial vector of degree at most \(d\).  Then \(h\) is SoS if
and only if it admits a Gram representation
$
    h(x)=z_d(x)^\top Q z_d(x)
     \text{ for some symmetric }Q\succeq0.
$
See, e.g.,
\cite{parrilo2000structured,lasserre2001global}.
\end{definition}

\begin{definition}[Nonconstant monomials and strict Gram representation]
\label{def:strict-gram}
We denote by
\[
    \textbf{v}_d(x)
    :=\bigl([x]_1^\top,\ldots,[x]_d^\top\bigr)^\top
\]
the vector of all nonconstant monomials of degrees \(1,\ldots,d\).  Its
dimension is
$
    N_{n,d}=\sum_{k=1}^d\binom{n+k-1}{k}
           =\binom{n+d}{d}-1.
$
When \(n\) and \(d\) are fixed, we write \(N:=N_{n,d}\).
For a polynomial \(q\) with zero constant and linear terms, we say that \(q\)
admits a \emph{strict Gram representation} on \(\textbf{v}_d\) if
\[
    q(x)=\textbf{v}_d(x)^\top Q\textbf{v}_d(x)
    \qquad\text{for some symmetric }Q\succ0.
\]
The matrix \(Q\) is called a \emph{Gram matrix} of \(q\).  
For a chosen strict Gram representation, \(\lambda_{\min}(Q)>0\) is called
its \emph{Gram margin}.  The Gram matrix need not be unique, so the margin is
attached to the chosen representation.

\end{definition}

\begin{definition}[SoS-convexity]
\label{def:sos-convex}
A polynomial \(h\) is \emph{SoS-convex} if
$
    y^\top\nabla^2 h(x)y
$
is an SoS polynomial jointly in \((x,y)\).  Equivalently, its Hessian is an
SoS matrix \cite{ahmadi2013complete}.  We use the standard implication that the first-order Bregman
remainder of an SoS-convex polynomial is SoS; the precise references are
recalled in \Cref{rem:standard-reference}.
\end{definition}

\section{How Regularization Creates SoS Exactness: An Overview}

\label{sec:proof-overview}

Before turning to the technical results, we give a short overview of the
mechanism behind the proof. The aim is to explain why increasing
regularization can change the SoS behaviour of the model even when the model
itself remains nonconvex.

The phenomenon is perhaps easiest to understand geometrically.  Suppose that
the Taylor model has negative curvature at the base point \(s=0\).  Since the
Hessian of the higher-order regularization term
\[
    \frac{\sigma}{2d}\|s\|^{2d},
    \qquad d\ge2,
\]
vanishes at the origin, increasing \(\sigma\) does not remove this negative
curvature.  Hence the regularized polynomial may remain locally nonconvex at
the base point for every finite value of \(\sigma\); see also
\cite{zhu2025global}.  Nevertheless, its algebraic behaviour can change
substantially as the regularization increases.  Interestingly,
Parrilo~\cite[Sec.~5.1]{parrilo2001minimizing} observed numerically that
certain quartically regularized polynomial models appear to be SoS once
the regularization is sufficiently large.

\begin{figure}[t]
    \centering
    \includegraphics[width=0.98\linewidth]{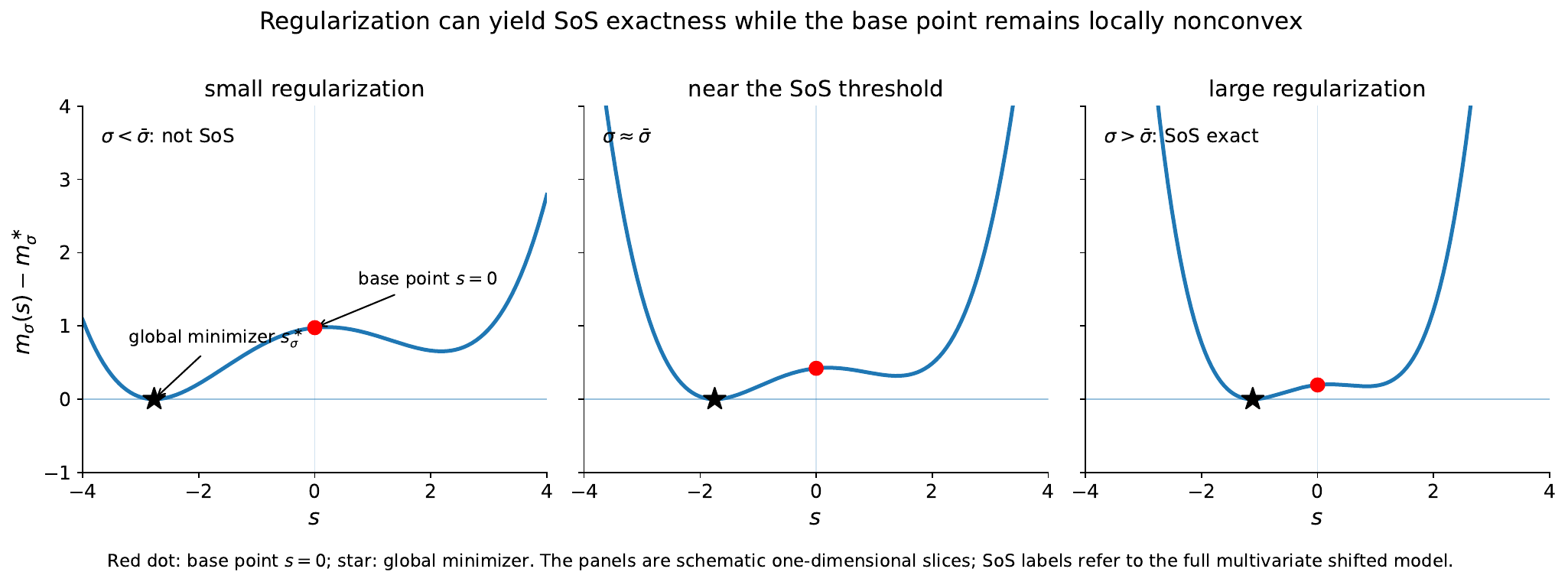}
    \caption{One-dimensional slices of multivariate
    regularized polynomial models as the regularization parameter increases.
    The red dot marks the base point $s=0$, where the negative quadratic
    curvature is unchanged by the higher-order regularizer; the star marks
    a global minimizer. }
    \label{fig:regularization-sos-transition}
\end{figure}

\Cref{fig:regularization-sos-transition} gives a
one-dimensional slice of a multivariate regularized polynomial at different
regularization strengths. The panels illustrate a change in geometry from ``not SoS'' to ``SoS exact'' as regularization increases. Note that the SoS labels refer to the
\emph{full multivariate optimality gap}, not to the plotted univariate
restriction: every nonnegative univariate polynomial is already SoS.

\paragraph{From polynomials to Taylor models.}
\label{par:coefficient-matching}

Any polynomial $P:\mathbb R^n\to\mathbb R$ of degree at most $p$
can be written exactly as a finite Taylor expansion,
\[
P(s)=\sum_{|\alpha|\le p}c_\alpha s^\alpha
=\sum_{j=0}^p\frac{1}{j!}T_j[s]^j,
\]
where $\alpha=(\alpha_1,\ldots,\alpha_n)$ is a multi-index,
$|\alpha|=\sum_{i=1}^n\alpha_i$,
$\alpha!=\prod_{i=1}^n\alpha_i!$, and
$s^\alpha=\prod_{i=1}^n s_i^{\alpha_i}$.
Here $T_j$ is the symmetric $j$th-order tensor determined by
the coefficients of $P$: for indices $i_1,\ldots,i_j$
with multiplicities $\alpha_1,\ldots,\alpha_n$, set
\[
(T_j)_{i_1\cdots i_j}=\alpha!\,c_\alpha.
\]
Indeed, each monomial $s^\alpha$ appears $j!/\alpha!$
times in the symmetric tensor contraction $T_j[s]^j$,
so the coefficient of $s^\alpha$ in
$\frac1{j!}T_j[s]^j$ is exactly $c_\alpha$.
In particular, $T_j=\nabla^jP(0)$; thus this representation
introduces no approximation or restriction on $P$. 

For convenience, throughout the paper we write the fixed
polynomial $P$ as its exact symmetric-tensor Taylor form
$\mathcal T_p$; the two symbols denote the same polynomial,
not two different models. 
We study the corresponding Euclidean-regularized polynomial
\begin{equation}
\label{eq:overview-regularized-model}
m_\sigma(s)
=
\mathcal T_{2d-1}(s)
+
\frac{\sigma}{2d}\|s\|^{2d},
\qquad d\ge2,
\end{equation}
where $\mathcal T_{2d-1}$ is any polynomial of degree at most
$2d-1$. This also includes even-degree Taylor models by
setting the $(2d-1)$st-order tensor to zero.
Writing
\[
\mathcal T_{2d-1}(s)
=
f_0+g^\top s+\frac12 H[s]^2
+\sum_{j=3}^{2d-1}\frac1{j!}T_j[s]^j,
\]
we have $T_0=f_0$, $T_1=g$, and $T_2=H$,
with $T_j[s]^j:=T_j[s,\ldots,s]$. 

\paragraph{Shifting the model at its global minimizer.}
The starting point of our proof is to rewrite the polynomial around a global
minimizer.  This is in the spirit of the classical optimality
characterizations used for trust-region and regularization subproblems;
see, e.g.,
\cite{wright2006numerical,cartis2022evaluation,zhu2025global}.
The difference here is that we use the resulting identity as an algebraic
decomposition for the SoS analysis. 
Let \(s_\sigma^*\) be a global minimizer of \(m_\sigma\), and define
\[
q_\sigma(v)
:=
m_\sigma(s_\sigma^*+v)-m_\sigma(s_\sigma^*).
\]
The polynomial \(q_\sigma\) is nonnegative and vanishes at \(v=0\).
The following identity is the starting point of the argument.

\begin{lemma}[Shifted optimality identity]
\label{lem:shifted-optimality-identity}
Let \(s_\sigma^*\) be a global minimizer of
\eqref{eq:intro-regularized-model}.  Then
\begin{equation}
\label{eq:overview-decomposition}
q_\sigma(v)
=
p_\sigma(v)
+
\frac{\sigma}{d}R_{s_\sigma^*,d}(v),
\end{equation}
where
\[
p_\sigma(v)
:=
\sum_{j=2}^{2d-1}
\frac{1}{j!}
\nabla^j\mathcal T_{2d-1}(s_\sigma^*)[v]^j \qquad
R_{s,d}(v)
:=
\frac12\|s+v\|^{2d}
-
\frac12\|s\|^{2d}
-
d\|s\|^{2d-2}s^\top v.
\]
\end{lemma}

\begin{proof}
Since the Taylor expansion of \(\mathcal T_{2d-1}\) around
\(s_\sigma^*\) is exact,
\[
q_\sigma(v)
=
\sum_{j=1}^{2d-1}
\frac{1}{j!}
\nabla^j\mathcal T_{2d-1}(s_\sigma^*)[v]^j
+
\frac{\sigma}{2d}
\left(
\|s_\sigma^*+v\|^{2d}
-
\|s_\sigma^*\|^{2d}
\right).
\]
The first-order optimality condition
\[
\nabla\mathcal T_{2d-1}(s_\sigma^*)
+
\sigma\|s_\sigma^*\|^{2d-2}s_\sigma^*
=
0
\]
cancels the linear term and gives
\eqref{eq:overview-decomposition}.
\end{proof}

We stress that \Cref{lem:shifted-optimality-identity} is an identity used in
the analysis; it does not assume that the global minimizer is known in
advance.  The purpose of the subsequent SoS representation is precisely to recover
this minimizer.

\paragraph{The three steps of the proof.}
The decomposition \eqref{eq:overview-decomposition}
separates the proof into two parts whose behaviour as \(\sigma\) increases is
very different.

\emph{First, the regularization term provides a positive SoS core.}
The polynomial \(R_{s,d}\) is precisely the first-order Bregman remainder of
$
\Phi_d(s)=\frac12\|s\|^{2d}.
$
Since \(\Phi_d\) is
SoS-convex, \(R_{s,d}\) is SoS. This is standard result of SoS-convexity
\cite{ahmadi2013complete,helton2010semidefinite}. However, the key result of
\Cref{sec The SoS property for the Bregman} is considerably stronger:
for every \(s\neq0\), this remainder admits a positive-definite Gram
representation on all nonconstant monomials up to degree \(d\).
Informally, the SoS representation therefore has some ``room'' to absorb
additional polynomial terms.  This is made precise through the strict Gram
margin developed in \Cref{sec The SoS property for the Bregman}.

\emph{Second, the remaining Taylor terms stay bounded.}
\(p_\sigma\) begins at quadratic order and
has degree at most \(2d-1\). Note that \(p_\sigma\) contains neither a constant nor a linear term; the constant term disappears because we shift by the value at
\(s_\sigma^*\), while the linear term is eliminated by the first-order
optimality condition.  
Its coefficients depend only on the fixed Taylor
tensors and on \(s_\sigma^*\).  As the regularization increases,
\(s_\sigma^*\) approaches the base point, and these coefficients remain
uniformly bounded; see \Cref{sec Main perturbation theorem}.

\emph{Finally, the regularization eventually dominates.}
Generally, as \(\sigma\) grows, the global minimizer moves towards the base point $s=0$. The
important question is how quickly this happens, since the Gram margin itself
also depends on \(s_\sigma^*\).  Under the assumption
$
g\neq0,
$ (i.e., the base point has not yet
reached first-order stationarity), we show that
\[
\|s_\sigma^*\|
=
\Theta\!\left(\sigma^{-1/(2d-1)}\right).   
\]
Hence the strict Gram margin contributed by
$
\frac{\sigma}{d}R_{s_\sigma^*,d}
$
grows with \(\sigma\), whereas the coefficients of \(p_\sigma\) remain
bounded.  For sufficiently large \(\sigma\), the regularization term therefore
absorbs the remaining Taylor terms and the entire shifted polynomial
\(q_\sigma\) becomes SoS.  The perturbation argument and the resulting SoS
exactness theorem are proved in \Cref{sec Main perturbation theorem}. For example, in the quartic case \(d=2\),
$
\|s_\sigma^*\|
=
\Theta(\sigma^{-1/3}), 
\sigma\|s_\sigma^*\|^2
=
\Theta(\sigma^{1/3}),
$
so the effective strict Gram margin grows like \(\sigma^{1/3}\). In one line, the proof structure can be summarized as
\[
\boxed{
\begin{gathered}
\text{shift at the minimizer}
\;\Longrightarrow\;
\text{bounded Taylor part + growing strict Gram part}\\
\Longrightarrow\;\text{SoS exactness as }\sigma\text{ increases.}
\end{gathered}}
\]

\section{\texorpdfstring{Bregman Remainders and Strict Gram Margins}{Bregman Remainders and Strict Gram Margins}}
\label{sec The SoS property for the Bregman}

This section studies the part of the model that comes from the
regularization term.  The reason is simple: after a regularized polynomial is
shifted around a stationary point, the contribution of the regularizer is its
first-order Bregman remainder.  Understanding this remainder will tell us how
much of the remaining Taylor polynomial the regularizer can absorb. 
Let \(d\ge2\) and define the Euclidean regularizer
\[
    \Phi_d(s):=\frac12\|s\|^{2d}.
\]
For a base point \(s\in\R^n\) and a displacement \(v\in\R^n\), its Bregman
remainder is
\begin{equation}
\label{eq:Rsd}
R_{s,d}(v)
:=
\Phi_d(s+v)-\Phi_d(s)-\nabla \Phi_d(s)^\top v
=
\frac12\|s+v\|^{2d}
-\frac12\|s\|^{2d}
-d\|s\|^{2d-2}s^\top v.
\end{equation}
It is standard that this polynomial is SoS when \(\Phi_d\) is SoS-convex;
see \Cref{rem:standard-reference}.  Our main goal is stronger: we show that
for every nonzero direction, the Bregman remainder has a \emph{strict} Gram
representation.  This extra room inside the SoS cone is what will later allow
us to absorb the other Taylor terms.

We use the monomial vector \(\textbf{v}_d\) and the notion of strict Gram
representation from \Cref{def:strict-gram}.

\begin{lemma}[SoS Bregman remainder]
\label{lem:bregman-SoS}
For every \(d\geq 2\), the polynomial \(\Phi_d(s)=\frac12\|s\|^{2d}\)
is SoS-convex. Therefore, its first-order
Bregman remainder $R_{s,d}(v)$ is a sum-of-squares polynomial. Specifically, for any $v \in \R^n$, let $w=s+v$,
\[
\Phi_d(w)-\Phi_d(s)-\nabla \Phi_d(s)^\top(w-s)
\]
is SoS jointly in \((s,w)\). This means that \(R_{s,d}(v)\) is SoS
jointly in \((s,v)\), and hence is \(R_{s,d}(v)\) SoS in \(v\) for every fixed \(s\).
\end{lemma}

\begin{proof}
\(\Phi_d\) is clearly
SoS-convex, to see this. We have
$
\nabla^2 \Phi_d(s)
=
d\|s\|^{2d-2}I
+
2d(d-1)\|s\|^{2d-4}ss^\top.
$
Therefore, for $w \in \R^n$
\[
w^\top\nabla^2 \Phi_d(s)w
=
d\|s\|^{2d-2}\|w\|^2
+
2d(d-1)\|s\|^{2d-4}(s^\top w)^2.
\]
Both terms are SoS polynomials in \((s,w)\), so 
\(\Phi_d\) is SoS-convex. The first-order characterization of SoS-convexity gives
\[
\Phi_d(w)-\Phi_d(s)-\nabla \Phi_d(s)^\top(w-s)\in\SoS.
\]
Equivalently, using the integral Taylor remainder,
\[
\Phi_d(w)-\Phi_d(s)-\nabla \Phi_d(s)^\top(w-s)
=
\int_0^1(1-t)(w-s)^\top
\nabla^2\Phi_d(s+t(w-s))(w-s)\,dt,
\]
and the integral remains SoS because the SoS cone is closed.
Setting \(w=s+v\) gives \eqref{eq:Rsd}.
\end{proof}

\begin{remark}[Standard reference]
\label{rem:standard-reference}
The result that 
\[
f \text{ is SoS convex i.e.,} \quad \nabla^2 f \text{ is an SoS matrix}
\quad\Longrightarrow\quad
f(y)-f(x)-\nabla f(x)^\top{(y-x)}\text{ is SoS}
\]
is standard. It appears as part of Theorem~3.1 in
Ahmadi and Parrilo~\cite{ahmadi2013complete}.
A closely related stationary-point statement appears as Lemma~8 in
Helton and Nie~\cite{helton2010semidefinite}: if \(f(u)=0\),
\(\nabla f(u)=0\), and \(\nabla^2 f\) is an SoS matrix, then \(f\) is SoS.
The complete statements of both results, in the form in which we
use them, are reproduced in \Cref{app:cited-sos-results}.  
We emphasize that these
cited results are used \emph{only} to establish \Cref{lem:bregman-SoS}. The main result of this paper and the perturbation argument in
\Cref{thm:eventual-SoS} requires a positive-definite strict Gram representation where \Cref{lem:bregman-SoS} is not enough. We will use the subsequent sections to develop the strict Gram argument.
\end{remark}

\subsection{\texorpdfstring{A Strict Gram Margin for Euclidean Regularization}{A Strict Gram Margin for Euclidean Regularization}}

Being SoS is not enough for the perturbation argument used later: the Gram
matrix could still be singular.  We now show that the Euclidean Bregman
remainder has a positive-definite Gram matrix whenever \(s\neq0\), and we
quantify how its smallest eigenvalue depends on \(\|s\|\).

\begin{theorem}
    [Strict Gram representation]
\label{lem:strict-gram}
For a fixed \(s\neq0\), $d \ge 2$ and $R_{s,d}(v)$ defined in \eqref{eq:Rsd}, there exists a symmetric matrix
\(Q_R(s)\succ0\) such that
\begin{equation}
\label{eq:gram-R}
R_{s,d}(v)=\textbf{v}_d(v)^\top Q_R(s)\textbf{v}_d(v).
\end{equation}
Moreover, there exists a constant \(c_{n,d}>0\), independent of
\(s\), such that for \(0<\|s\|\leq1\),
\begin{equation}
\label{eq:margin}
\lmin(Q_R(s))
\geq
c_{n,d}\|s\|^{2d-2}.
\end{equation}
\end{theorem}

\begin{proof}[Proof of \Cref{lem:strict-gram}]
Applying the integral Taylor formula to \eqref{eq:Rsd} gives
\begin{align}
R_{s,d}(v)
=
 d\underbrace{\int_0^1(1-t)\|s+tv\|^{2d-2}\|v\|^2\,dt}_{I_1(v)}
+2d(d-1)
\underbrace{\int_0^1(1-t)\|s+tv\|^{2d-4}
\bigl((s+tv)^\top v\bigr)^2\,dt}_{I_2(v)}.
\label{rsd integral}
\end{align}

Here \(I_1\) and \(I_2\) denote the two integrals displayed in
\eqref{rsd integral}.  Their roles are different: \(I_1\) supplies the
strict part of the Gram matrix, whereas \(I_2\) contributes only a
positive-semidefinite part.  The remaining lemmas verify these two facts and
then transport the result from \(e_1\) to an arbitrary nonzero point.

\begin{center}
\begin{tabular}{c|c|c}
term & role in the proof & reference\\ \hline
\(I_1\) & positive-definite Gram contribution &
\Cref{lem:first-integral-strict-gram,lem:radial-feature-spanning}\\
\(I_2\) & positive-semidefinite Gram contribution &
\Cref{lem:second-integral-psd}
\end{tabular}
\end{center}

By \Cref{lem:radial-feature-spanning}, the polynomials in the first integral at
\(s=e_1\) have the full coefficient-spanning property. Hence
\Cref{lem:first-integral-strict-gram} gives a Gram matrix
\(Q_1(e_1)\succ0\) satisfying
\[
I_1(v)=\textbf{v}_d(v)^\top Q_1(e_1)\textbf{v}_d(v).
\]
By \Cref{lem:second-integral-psd}, the second integral has a Gram matrix
\(Q_2(e_1)\succeq0\) satisfying
\[
I_2(v)=\textbf{v}_d(v)^\top Q_2(e_1)\textbf{v}_d(v).
\]
Therefore, using \eqref{rsd integral}, the Gram matrix for \(R_{e_1,d}\) satisfies
\begin{equation}
\label{eq:assembled-gram-at-e1}
Q_R(e_1) 
:= 
dQ_1(e_1)+2d(d-1)Q_2(e_1) 
\succeq d\lambda_{\min}(Q_1(e_1)) I_N \succ 0. 
\end{equation} 
Applying \Cref{lem:gram-transport-scaling} gives the change of coordinates from $e_1$ to any fixed \(s\neq0\),  
\[
\lambda_{\min}(Q_R(s))
\geq
\|s\|^{2d-2}\lambda_{\min}(Q_R(e_1)),
\qquad
0<\|s\|\leq1.
\]
Thus \eqref{eq:margin} holds with
$
c_{n,d}:=\lambda_{\min}(Q_R(e_1))>0.
$
\end{proof}

 \subsubsection{\texorpdfstring{The Two Integral Terms and the Spanning Argument}{The Two Integral Terms and the Spanning Argument}}
\label{sec tech lemmas}

We first dispose of the easier term: \(I_2\) is a weighted sum of squares and
therefore contributes a positive-semidefinite Gram matrix.

\begin{lemma}[The second integral contribution is SoS]
\label{lem:second-integral-psd}
For every fixed \(s\) and \(d\ge2\), there exists a symmetric matrix
\(Q_2(s)\succeq0\) such that
\begin{equation}
\label{eq:second-contribution-gram}
I_2(v)=\textbf{v}_d(v)^\top Q_2(s)\textbf{v}_d(v).
\end{equation}
\end{lemma}

\begin{proof}
For \(d\ge2\), the multinomial expansion gives
\begin{equation}
\label{eq:second-integrand-squares}
\|s+tv\|^{2d-4}\bigl((s+tv)^\top v\bigr)^2
=
\sum_{|\alpha|=d-2}\widetilde c_\alpha
\bigl[\prod_{j=1}^n(s_j+tv_j)^{\alpha_j} (s+tv)^\top v\bigr]^2,
\end{equation}
where
\[
\widetilde c_\alpha
=
\binom{d-2}{\alpha}
=
\frac{(d-2)!}{\alpha_1!\cdots\alpha_n!}>0.
\]
For each \(\alpha\), let \(\widetilde a_\alpha(s,t)\)
denote the coefficient vector of
$
\prod_{j=1}^n(s_j+tv_j)^{\alpha_j} (s+tv)^\top v
$
with respect to \(\textbf{v}_d(v)\).  Define
\begin{equation}
\label{eq:q2-definition}
Q_2(s)
:=
\sum_{|\alpha|=d-2}\widetilde c_\alpha
\int_0^1(1-t)\,
\widetilde a_\alpha(s,t)\widetilde a_\alpha(s,t)^\top\,dt.
\end{equation}
Then \eqref{eq:second-contribution-gram} follows from
\eqref{eq:second-integrand-squares}.  Moreover, for every
\(\xi\in\R^N\),
\begin{equation}
\label{eq:q2-psd-check}
\xi^\top Q_2(s)\xi
=
\sum_{|\alpha|=d-2}\widetilde c_\alpha
\int_0^1(1-t)
\bigl(\xi^\top\widetilde a_\alpha(s,t)\bigr)^2\,dt
\ge0.
\end{equation}
Thus \(Q_2(s)\succeq0\).  In particular, each matrix under the integral is
rank at most one and positive semidefinite, and the strictly positive
coefficients \(2d(d-1)\), \(\widetilde c_\alpha\), and \(1-t\) prevent any
cancellation between the square contributions.
\end{proof}

The first integral \(I_1\) is the source of strictness.  To prove this, we
track the coefficients of the squared polynomials appearing in the integral
and ask whether they cover every nonconstant monomial direction.

Throughout the following lemmas, coefficient vectors are taken with respect
to the fixed ordering of \(\textbf{v}_d\) from \Cref{def:strict-gram}.

\begin{lemma}[The first integral contribution and coefficient spanning]
\label{lem:first-integral-strict-gram}
The following statements hold.
\begin{enumerate}
\item For fixed \(s\), let
$
\mathfrak{p}_{\alpha,i,t}(v):=v_i\prod_{j=1}^n(s_j+tv_j)^{\alpha_j}.
$
The first integral in \eqref{rsd integral} can be expressed as
\begin{equation}
\label{eq:first-contribution-gram}
I_1(v)
=
\sum_{i=1}^n\sum_{|\alpha|=d-1}c_\alpha
\int_0^1(1-t)\mathfrak{p}_{\alpha,i,t}(v)^2\,dt
=
\textbf{v}_d(v)^\top Q_1(s)\textbf{v}_d(v),
\end{equation}
where
\begin{equation}
\label{eq:q1-definition}
Q_1(s)
:=
\sum_{i=1}^n\sum_{|\alpha|=d-1}c_\alpha
\int_0^1(1-t)\,
a_{i,\alpha}(s,t)a_{i,\alpha}(s,t)^\top\,dt.
\end{equation}

\item For fixed \(s\in\R^n\), let \(a_{i,\alpha}(s,t)\in\R^N\) denote the
coefficient vector of \(\mathfrak{p}_{\alpha,i,t}\) with respect to the
ordered monomial vector \(\textbf{v}_d(v)\); that is,
\(\mathfrak{p}_{\alpha,i,t}(v)=a_{i,\alpha}(s,t)^\top\textbf{v}_d(v)\).
Equivalently, we write
\[
a_{i,\alpha}(s,t)
=
\coeff_{\textbf{v}_d}\bigl(\mathfrak{p}_{\alpha,i,t}\bigr)
\in\R^N.
\]
The map \(t\mapsto a_{i,\alpha}(s,t)\) is a vector-valued polynomial\footnote{For every fixed \(t\), \(a_{i,\alpha}(s,t)\) is a vector in
\(\R^N\).} of
degree at most \(d-1\) in \(t\). In particular,
\begin{eqnarray}
  \mathfrak{p}_{\alpha,i,t}(v)
=
a_{i,\alpha}(s,t)^\top\textbf{v}_d(v)
=
\left(
\sum_{k=0}^{d-1}b_{i,\alpha,k}(s)t^k
\right)^\top\textbf{v}_d(v).
\label{coeff vector b}  
\end{eqnarray}
where for fixed $s$, each $b_{i,\alpha,k} \in \R^N$. 

\item \(Q_1(s)\succ0\) if and only if $\{b_{i,\alpha,k}(s)\}_{i,\alpha,k}$ satisfies the following condition:
\begin{equation}
\label{eq:annihilator-spanning-condition}
\forall 
 0\neq\xi\in\R^N , \quad
\exists\ i=1,\ldots,n,\quad |\alpha|=d-1,\quad 0\leq k\leq d-1
\quad\text{s.t.}\quad
\xi^\top b_{i,\alpha,k}(s)\neq0.
\tag{Spanning property}
\end{equation}
\end{enumerate}
\end{lemma}

\begin{proof}
Expanding term by term gives
\begin{align*}
\|s+tv\|^{2d-2}\|v\|^2
 &=
\bigg(\sum_{j=1}^n(s_j+tv_j)^2\bigg)^{d-1}
\bigg(\sum_{i=1}^n v_i^2\bigg)  =
\bigg(
\sum_{|\alpha|=d-1}c_\alpha
\prod_{j=1}^n(s_j+tv_j)^{2\alpha_j}
\bigg)
\bigg(\sum_{i=1}^n v_i^2\bigg) \\
 &=
\sum_{i=1}^n\sum_{|\alpha|=d-1}c_\alpha
\bigg[
v_i\prod_{j=1}^n(s_j+tv_j)^{\alpha_j}
\bigg]^2,
\quad \text{where} \quad
c_\alpha
=
\binom{d-1}{\alpha}
=
\frac{(d-1)!}{\alpha_1!\cdots\alpha_n!}>0.
\end{align*}
Therefore,
\begin{equation}
\label{eq:first-contribution-squares}
I_1(v)
=
\sum_{i=1}^n\sum_{|\alpha|=d-1}c_\alpha
\int_0^1(1-t)\mathfrak{p}_{\alpha,i,t}(v)^2\,dt.
\end{equation}
This proves the first equality in \eqref{eq:first-contribution-gram}. By the
definition of \(a_{i,\alpha}(s,t)\),
\[
\mathfrak{p}_{\alpha,i,t}(v)^2
=
\textbf{v}_d(v)^\top
a_{i,\alpha}(s,t)a_{i,\alpha}(s,t)^\top
\textbf{v}_d(v),
\]
and hence \eqref{eq:q1-definition} proves the second equality in
\eqref{eq:first-contribution-gram}. This proves part~1. Part~2 follows because,
for fixed \(s\), every entry of \(a_{i,\alpha}(s,t)\) is a polynomial in
\(t\) of degree at most \(d-1\). 
For any \(\xi\in\R^N\), we have
\begin{equation}
\label{eq:q1-quadratic-form}
\xi^\top Q_1(s)\xi
=
\sum_{i=1}^n\sum_{|\alpha|=d-1}c_\alpha
\int_0^1(1-t)
\bigl[\xi^\top a_{i,\alpha}(s,t)\bigr]^2\,dt.
\end{equation}
Since every \(c_\alpha>0\), the weight \(1-t\) is strictly positive on
\((0,1)\), and every integrand in \eqref{eq:q1-quadratic-form} is
nonnegative, \(Q_1(s)\succeq0\).

We now prove part 3. 
Suppose first that \eqref{eq:annihilator-spanning-condition} holds, and fix any
\(0\neq\xi\in\R^N\). Then there exist \(i,\alpha,k\) such that
\(\xi^\top b_{i,\alpha,k}(s)\neq0\). $ \xi^\top a_{i,\alpha}(s,t)$ is a scalar polynomial with some non-zero coefficient (\(\xi^\top b_{i,\alpha,k}(s)\neq0\)) and is therefore not identically zero
\[
 \xi^\top a_{i,\alpha}(s,t)
=
\sum_{\ell=0}^{d-1}\xi^\top b_{i,\alpha,\ell}(s)t^\ell \neq 0.
\]
Therefore, for any fixed $\xi \in \R^N$,  the following integrant is a nonnegative (not identically zero) polynomial in $t \in (0, 1)$. Thus, we have
$$
\int_0^1(1-t)
\bigl[\xi^\top a_{i,\alpha}(s,t)\bigr]^2\,dt>0.
$$ 
It follows from \eqref{eq:q1-quadratic-form} that
\(\xi^\top Q_1(s)\xi>0\). Since this holds for every nonzero \(\xi\),
\(Q_1(s)\succ0\).

Conversely, suppose that \eqref{eq:annihilator-spanning-condition} fails.
Then there exists \(0\neq\xi\in\R^N\) such that
\[
\xi^\top b_{i,\alpha,k}(s)=0
\qquad
\text{for every }i,\alpha,k.
\]
This gives
\(\xi^\top a_{i,\alpha}(s,t)=0\) for every \(i,\alpha\) and every
\(t\in(0,1)\). Hence \eqref{eq:q1-quadratic-form} yields
\(\xi^\top Q_1(s)\xi=0\), so \(Q_1(s)\) is not positive definite. This
proves part~3.

\end{proof}

The next lemma verifies the spanning property for Euclidean
regularization. The low-degree
examples in \Cref{example small case,example small case2} give illustrations of this lemma.

\begin{lemma}[Spanning for the Euclidean radial polynomials]
\label{lem:radial-feature-spanning}
Fix $s=re_1$ with $r>0$.  The coefficient vectors $\{b_{i,\alpha,k}(re_1)\}_{i,\alpha,k}$ are defined by \eqref{coeff vector b}: 
$$
a_{i,\alpha}(re_1,t)^\top\textbf{v}_d(v)
=
\left(
\sum_{k=0}^{d-1}b_{i,\alpha,k}(re_1)t^k
\right)^\top\textbf{v}_d(v). 
$$
Then $\{b_{i,\alpha,k}(re_1)\}_{i,\alpha,k}$ satisfies \eqref{eq:annihilator-spanning-condition}; consequently, $Q_1(re_1)\succ0$. 
\end{lemma}

\begin{proof}   
Fix any monomial \(v^\beta\) with
$
1\leq|\beta|\leq d.
$
Choose an index \(i\) such that \(\beta_i>0\), and set
$
\gamma:=\beta-e_i.
$
Thus \(|\gamma|=|\beta|-1\leq d-1\). Define a multi-index \(\alpha\) by
\[
\alpha_j=\gamma_j\quad (j=2,\ldots,n),
\qquad
\alpha_1=d-1-\sum_{j=2}^n\gamma_j.
\]
Then \(|\alpha|=d-1\), and moreover
$\alpha_1-\gamma_1
=
d-1-|\gamma|
=
d-|\beta|
\geq0.
$
Expanding the corresponding polynomial gives
\begin{align}
\mathfrak{p}_{\alpha,i,t}(v)
&: =
v_i(r+tv_1)^{\alpha_1}
\prod_{j=2}^n(tv_j)^{\gamma_j}  =
\sum_{k=0}^{\alpha_1}
\binom{\alpha_1}{k}
r^{\alpha_1-k}
t^{k+\sum_{j=2}^n\gamma_j}
v_i v_1^k\prod_{j=2}^n v_j^{\gamma_j}.
\label{eq:radial-feature-expansion}
\end{align}
The coefficient of \(t^{|\gamma|}=t^{|\beta|-1}\) is obtained by taking
\(k=\gamma_1\), and equals
\begin{equation}
\label{eq:target-monomial-coefficient}
\binom{\alpha_1}{\gamma_1}
r^{\alpha_1-\gamma_1}
v_i v_1^{\gamma_1}\prod_{j=2}^n v_j^{\gamma_j}
=
\binom{\alpha_1}{\gamma_1}
r^{\alpha_1-\gamma_1}
v^\beta.
\end{equation}

Comparing the coefficient of \(t^{|\beta|-1}\) in
\eqref{eq:radial-feature-expansion} with
\eqref{coeff vector b} gives
\begin{equation}
\label{eq:target-coefficient-vector}
b_{i,\alpha,|\beta|-1}(re_1)
=
\binom{\alpha_1}{\gamma_1}
r^{d-|\beta|}
\coeff_{\textbf{v}_d}(v^\beta).
\end{equation}
Now let \(0\neq \xi\in\R^N\). Since the coefficient vectors
$
    \left\{
    \coeff_{\textbf{v}_d}(v^\beta):
    1\le |\beta|\le d
    \right\}
$
form a basis of \(\R^N\), there always exists a monomial \(v^\beta\) such that
\[
    \xi^\top\coeff_{\textbf{v}_d}(v^\beta)\neq0.
\]
Choose \(i\) such that \(\beta_i>0\), and  
$
    \gamma:=\beta-e_i, 
    \alpha:=\gamma+(d-|\beta|)e_1.
$
Using \eqref{eq:target-coefficient-vector} and since 
$
    \binom{\alpha_1}{\gamma_1}r^{d-|\beta|}>0,
$
we obtain
\[
    \xi^\top b_{i,\alpha,|\beta|-1}(re_1)
    =
    \binom{\alpha_1}{\gamma_1}
    r^{d-|\beta|}
    \xi^\top\coeff_{\textbf{v}_d}(v^\beta)
    \neq0.
\]
Therefore
\(\{b_{i,\alpha,k}(re_1)\}_{i,\alpha,k}\) satisfies
\eqref{eq:annihilator-spanning-condition}.  Consequently, by \Cref{lem:first-integral-strict-gram},
$
    Q_1(re_1)\succ0.
$

\end{proof}

It remains only to move the certificate from \(e_1\) to
an arbitrary nonzero \(s\).  Orthogonal invariance handles the direction and
homogeneity handles the length.

\begin{lemma}[Orthogonal transport and homogeneous scaling]
\label{lem:gram-transport-scaling}
According to \Cref{lem:radial-feature-spanning}, 
\[
R_{e_1,d}(v)=\textbf{v}_d(v)^\top Q(e_1)\textbf{v}_d(v),
\qquad
Q(e_1)\succ0.
\]
Then every \(s\neq0\) admits a positive-definite Gram matrix \(Q_R(s)\).
Moreover, if \(s=ru\), \(\|u\|=1\), and \(0<r\le1\), then
\begin{equation}
\label{eq:transport-scaling-bound}
\lambda_{\min}(Q_R(ru))
\ge
r^{2d-2}\lambda_{\min}(Q(e_1)).
\end{equation}
\end{lemma}

%
%
\begin{proof}
First transport the certificate from \(e_1\) to an arbitrary unit direction.
For \(u\in\mathbb S^{n-1}\), choose an orthogonal matrix \(O_u\) such that
$
u=O_u e_1.
$
The Euclidean radial structure gives
$
R_{u,d}(v)=R_{e_1,d}(O_u^\top v).
$
With the orthogonally normalized monomial blocks in
\Cref{def:strict-gram}, there exists a block-orthogonal matrix
\(T_{O_u}\) satisfying
\begin{equation}
\label{eq:orthogonal-feature-action}
\textbf{v}_d(O_u^\top v)=T_{O_u}\textbf{v}_d(v),
\qquad
T_{O_u}^\top T_{O_u}=I.
\end{equation}
Define $
Q_R(u):=T_{O_u}^\top Q(e_1)T_{O_u}.
$
Then \(Q_R(u)\) is a Gram matrix for \(R_{u,d}\), and
\begin{equation}
\label{eq:unit-sphere-margin}
\lambda_{\min}(Q_R(u))
=
\lambda_{\min}(Q(e_1))
=:\mu_{n,d}>0.
\end{equation}
Thus the orthogonal matrix \(O_u\) transports the certificate between unit
directions and preserves its Gram margin.
Now let \(s=ru\), where \(r>0\).  By homogeneity,
\begin{equation}
\label{eq:remainder-homogeneous-scaling}
R_{ru,d}(v)
=
r^{2d}R_{u,d}\left(\frac vr\right).
\end{equation}
Since the degree-\(k\) homogeneous block of \(\textbf{v}_d(v/r)\) equals
\(r^{-k}[v]_k\), define the block-diagonal matrix
\begin{equation}
\label{eq:radial-block-scaling}
E_r
:=
\operatorname{diag}
\left(
 r^{d-1}I_{[v]_1},
 r^{d-2}I_{[v]_2},
 \ldots,
 rI_{[v]_{d-1}},
 I_{[v]_d}
\right),
\qquad
r^d\textbf{v}_d\left(\frac vr\right)=E_r\textbf{v}_d(v).
\end{equation}
It follows that
\begin{equation}
\label{eq:scaled-gram-definition}
R_{ru,d}(v)
=
\textbf{v}_d(v)^\top E_rQ_R(u)E_r\textbf{v}_d(v),
\qquad
Q_R(ru):=E_rQ_R(u)E_r.
\end{equation}
Since \(E_r\) is invertible, \(Q_R(ru)\succ0\) for every \(r>0\).  If
\(0<r\le1\), then \(\lambda_{\min}(E_r)=r^{d-1}\), and therefore
\[
\lambda_{\min}(Q_R(ru))
\ge
\lambda_{\min}(E_r)^2\lambda_{\min}(Q_R(u))
=
r^{2d-2}\mu_{n,d}.
\]
This proves \eqref{eq:transport-scaling-bound}.
\end{proof}

{

\begin{example}[Example proof $Q_1 \succ 0$ for \(d=2,n=3\)]

\label{example small case}
For \(d=2\), we have
\[
\|s+tv\|^{2}\|v\|^2
=
\sum_{i=1}^3\sum_{j=1}^3
\bigl[v_i(s_j+t v_j)\bigr]^2,
\]
and hence
\[
{2I_1(v)}
=
2\sum_{i,j=1}^3
\int_0^1(1-t)\bigl[v_i(s_j+t v_j)\bigr]^2\,dt .
\]
For each \(i,j\), let \(a_{ij}(s,t)\) be the coefficient vector of
\(v_i(s_j+t v_j)\) with respect to the nonconstant monomial vector
\[
z_2(v)
=
\begin{pmatrix}
v_1&v_2&v_3&
v_1^2&v_1v_2&v_1v_3&
v_2^2&v_2v_3&v_3^2
\end{pmatrix}^{\top}.
\]
Then
\[
Q_1(s)
=
2\sum_{i=1}^3\sum_{j=1}^3
\int_0^1(1-t)\,
a_{ij}(s,t)a_{ij}(s,t)^\top\,dt,
\qquad
{2I_1(v)=z_2(v)^\top Q_1(s)z_2(v)}.
\]
It suffices to consider \(s=re_1=(r,0,0)\), \(r>0\). Then
\[
v_i(s_j+t v_j)
=
\begin{cases}
v_i(r+t v_1), & j=1,\\
t v_i v_j, & j=2,3.
\end{cases}
\]
With respect to the above ordering of \(z_2(v)\), we have
\[
\begin{array}{lll}
a_{11}(t)=r e_{v_1}+t e_{v_1^2},&
a_{21}(t)=r e_{v_2}+t e_{v_1v_2},&
a_{31}(t)=r e_{v_3}+t e_{v_1v_3},\\[1mm]
a_{12}(t)=t e_{v_1v_2},&
a_{22}(t)=t e_{v_2^2},&
a_{32}(t)=t e_{v_2v_3},\\[1mm]
a_{13}(t)=t e_{v_1v_3},&
a_{23}(t)=t e_{v_2v_3},&
a_{33}(t)=t e_{v_3^2}.
\end{array}
\]
Note that since $\alpha = \{e_1, e_2, e_3\}$ here, we denote $a_{i,\alpha}(re_1,t)  = a_{ij}(t)$ for abbreviation. 
With respect to the above ordering of \(z_2(v)\), the nine coefficient
vectors $\{a_{11}(t), a_{21}(t), a_{31}(t) \dotsc a_{33}(t)\}$ are
\[
\begin{array}{lllllllll}
\begin{pmatrix} r\\0\\0\\ t\\0\\0\\0\\0\\0 \end{pmatrix},
&
\begin{pmatrix} 0\\r\\0\\ 0\\t\\0\\0\\0\\0 \end{pmatrix},
&
\begin{pmatrix} 0\\0\\r\\ 0\\0\\t\\0\\0\\0 \end{pmatrix},
&
\begin{pmatrix} 0\\0\\0\\ 0\\t\\0\\0\\0\\0 \end{pmatrix},
&
\begin{pmatrix} 0\\0\\0\\ 0\\0\\0\\t\\0\\0 \end{pmatrix},
&
\begin{pmatrix} 0\\0\\0\\ 0\\0\\0\\0\\t\\0 \end{pmatrix},
&
\begin{pmatrix} 0\\0\\0\\ 0\\0\\t\\0\\0\\0 \end{pmatrix},
&
\begin{pmatrix} 0\\0\\0\\ 0\\0\\0\\0\\t\\0 \end{pmatrix},
&
\begin{pmatrix} 0\\0\\0\\ 0\\0\\0\\0\\0\\t \end{pmatrix}.
\end{array}
\]
Since \(t\in(0,1)\), the last six vectors with \(j=2,3\), 
directly generate the five distinct quadratic directions
\[
e_{v_1v_2},\quad e_{v_1v_3},\quad
e_{v_2^2},\quad e_{v_2v_3},\quad e_{v_3^2},
\]
where the direction \(e_{v_2v_3}\) appears twice, through both
\(a_{32}(t) =a_{23}(t)=t e_{v_2v_3}\).
The remaining quadratic direction \(e_{v_1^2}\) is obtained from the
\(j=1\) vectors. Indeed, choosing two distinct values
\(t_0,t_1\in(0,1)\), we have
\[
a_{11}(t_1)-a_{11}(t_0)
=
(t_1-t_0)e_{v_1^2}.
\]
Hence \(e_{v_1^2}\) also belongs to the span. Therefore all six quadratic
directions
\[
e_{v_1^2},\quad e_{v_1v_2},\quad e_{v_1v_3},\quad
e_{v_2^2},\quad e_{v_2v_3},\quad e_{v_3^2}
\]
belong to the span.  Finally, since
\[
a_{11}(t_0)=r e_{v_1}+t_0e_{v_1^2},\qquad
a_{21}(t_0)=r e_{v_2}+t_0e_{v_1v_2},\qquad
a_{31}(t_0)=r e_{v_3}+t_0e_{v_1v_3},
\]
and \(r>0\), the linear directions
\[
e_{v_1},\quad e_{v_2},\quad e_{v_3}
\]
also belong to the span. Hence, the vectors \(a_{ij}(re_1,t)\), \(i,j=1,2,3\), \(t\in(0,1)\), span the full coefficient space generated by monomials of degrees \(1\) and \(2\). Namely,
\[
\boxed{
\operatorname{span}\{a_{ij}(re_1,t):\, i,j=1,2,3,\ t\in(0,1)\}
=
\operatorname{span}\{e_{v_1},e_{v_2},e_{v_3},
e_{v_1^2},e_{v_1v_2},e_{v_1v_3},
e_{v_2^2},e_{v_2v_3},e_{v_3^2}\}
=
\mathbb{R}^9 .
}
\]
Therefore, for any coefficient vector \(\xi\),
\[
\xi^\top Q_1(s)\xi
=
2\sum_{i=1}^3\sum_{j=1}^3
\int_0^1(1-t)
\bigl(\xi^\top a_{ij}(s,t)\bigr)^2\,dt .
\]
If \(\xi^\top Q_1(s)\xi=0\), then all nonnegative integrands vanish
identically, so \(\xi\) is orthogonal to all vectors \(a_{ij}(s,t)\).
Since these vectors span the full coefficient space, this implies
\(\xi=0\). Thus
$
Q_1(s)\succ0 .
$

\end{example}

\begin{example}[\(d=3,n=3\)]
\label{example small case2}
Let \(s=re_1=(r,0,0)\), \(r>0\).  For \(d=3\), we have
\(|\alpha|=d-1=2\).  We use the nonconstant monomial vector
\[
\begin{aligned}
z_3(v)
=
\big(
v_1,v_2,v_3,,
v_1^2,v_1v_2,v_1v_3,v_2^2,v_2v_3,v_3^2, 
v_1^3,v_1^2v_2,v_1^2v_3,
v_1v_2^2,v_1v_2v_3,v_1v_3^2,
v_2^3,v_2^2v_3,v_2v_3^2,v_3^3
\big)^\top .
\end{aligned}
\]
For each \(i=1,2,3\) and each multi-index \(|\alpha|=2\), let
\(a_{i,\alpha}(t)\) denote the coefficient vector of
\[
    v_i(re_1+tv)^\alpha
\]
with respect to \(z_3(v)\).  The possible multi-indices are
\[
\alpha : = \{
2e_1,\quad e_1+e_2,\quad e_1+e_3,\quad
2e_2,\quad e_2+e_3,\quad 2e_3 \}.
\]
The corresponding polynomials and coefficient vectors have the form
\[
\begin{aligned}
v_i(re_1+tv)^{2e_1}
&=
r^2v_i+2rt\,v_iv_1+t^2v_iv_1^2,
&
a_{i,2e_1}(t)
&=
r^2 e_{v_i}
+
2rt\, e_{v_iv_1}
+
t^2 e_{v_iv_1^2},
\\
v_i(re_1+tv)^{e_1+e_2}
&=
rt\,v_iv_2+t^2v_iv_1v_2,
&
a_{i,e_1+e_2}(t)
&=
rt\, e_{v_iv_2}
+
t^2 e_{v_iv_1v_2},
\\
v_i(re_1+tv)^{e_1+e_3}
&=
rt\,v_iv_3+t^2v_iv_1v_3,
&
a_{i,e_1+e_3}(t)
&=
rt\, e_{v_iv_3}
+
t^2 e_{v_iv_1v_3},
\\
v_i(re_1+tv)^{2e_2}
&=
t^2v_iv_2^2,
&
a_{i,2e_2}(t)
&=
t^2 e_{v_iv_2^2},
\\
v_i(re_1+tv)^{e_2+e_3}
&=
t^2v_iv_2v_3,
&
a_{i,e_2+e_3}(t)
&=
t^2 e_{v_iv_2v_3},
\\
v_i(re_1+tv)^{2e_3}
&=
t^2v_iv_3^2,
&
a_{i,2e_3}(t)
&=
t^2 e_{v_iv_3^2}.
\end{aligned}
\]
Note that we denote $a_{i,\alpha}(re_1,t)  = a_{i,\alpha}(t)$ for abbreviation. 
Thus the \(t^0\)-coefficients generate the degree-one directions
\[
e_{v_1},\quad e_{v_2},\quad e_{v_3}.
\]

More explicitly, each vector \(a_{i,\alpha}(t)\) is a vector-valued polynomial
of degree at most two:
\[
    a_{i,\alpha}(t)=b_0+t b_1+t^2 b_2 .
\]
Let \(\xi\in\R^{19}\) be orthogonal to
\(a_{i,\alpha}(t)\) for every \(i,\alpha\) and every
\(t\in(0,1)\).  Then
\[
 t\longmapsto \xi^\top a_{i,\alpha}(t)
\]
is a scalar polynomial which vanishes on \((0,1)\), and hence it vanishes
identically.  Therefore \(\xi\) is orthogonal to every coefficient vector
\(b_0,b_1,b_2\).  By the displayed expansions above, the \(t^0\)-,
\(t^1\)-, and \(t^2\)-coefficients generate, respectively, all monomial
coefficient directions of degrees one, two, and three.  Thus \(\xi=0\),
which proves that the value vectors \(a_{i,\alpha}(t)\), with
\(t\in(0,1)\), span the full coefficient space.

Therefore,
\[
\boxed{
\operatorname{span}\{a_{i,\alpha}(t):
i=1,2,3,\ |\alpha|=2,\ t\in(0,1)\}
=
\operatorname{span}\{e_{v^\beta}:1\le |\beta|\le 3\} = \R^{19}.
}
\]
That is, the vectors \(a_{i,\alpha}(t)\) span the full nonconstant coefficient
space generated by monomials of degrees \(1,2\), and \(3\).
\end{example}
}


\paragraph{The radial term $I_2$ is PSD but has no strict Gram margin}
Although $I_2(v)$ in \eqref{rsd integral} contains strong radial terms, it is not, by itself, 
\textbf{strictly positive} on the full nonconstant monomial vector \(\textbf{v}_d(v)\). 
For reference, when \(d=2\),
\[
I_1(v)
=
\frac12 \|s\|^2\|v\|^2
+\frac13 (s^\top v)\|v\|^2
+\frac1{12}\|v\|^4, \qquad
I_2(v)
=
\frac12 (s^\top v)^2
+\frac13 (s^\top v)\|v\|^2
+\frac1{12}\|v\|^4.
\]
When \(d=3\),
\[
\begin{aligned}
I_1(v)
=
\frac12 \|s\|^4\|v\|^2
+\frac23 \|s\|^2(s^\top v)\|v\|^2
+\frac16 \|s\|^2\|v\|^4  +
+\frac13 (s^\top v)^2\|v\|^2
+\frac15 (s^\top v)\|v\|^4
+\frac1{30}\|v\|^6,
\end{aligned}
\]
whereas
\[
\begin{aligned}
I_2(v)
=
\frac12 \|s\|^2(s^\top v)^2
+\frac13 \|s\|^2(s^\top v)\|v\|^2
+\frac1{12}\|s\|^2\|v\|^4 +
+\frac13 (s^\top v)^3
+\frac5{12}(s^\top v)^2\|v\|^2
+\frac15 (s^\top v)\|v\|^4
+\frac1{30}\|v\|^6 .
\end{aligned}
\]
Thus \(I_2\) indeed contributes more visibly radial terms, such as
\[
(s^\top v)^2
\qquad\text{and}\qquad
\|s\|^2(s^\top v)^2 .
\]
By contrast, \(I_1\) contains the isotropic factor
\(\|v\|^2\) and its associated squared terms span all nonconstant monomial
directions. Hence \(I_1\), rather than \(I_2\), is the source of the strict
Gram margin. More generally, we can prove that if \(n\ge2\), take \(w\perp s\), \(\|w\|=1\), and \(v=\tau w\).
Then
$
(s+t\tau w)^\top(\tau w)=t\tau^2,
$and hence
\[
I_2(\tau w)=O(\tau^4)
\qquad\text{as }\tau\to0.
\]
On the other hand, since \(\textbf{v}_d(\tau w)\) contains linear monomials,
\[
\|\textbf{v}_d(\tau w)\|^2=O(\tau^2).
\]
Thus no bound of the form
\[
I_2(v)\ge \mu \|\textbf{v}_d(v)\|^2,
\qquad \mu>0,
\]
can hold. Therefore \(I_2\) is positive semidefinite but not strictly
positive on the full coefficient space. The strict Gram margin is instead
provided by \(I_1\).

More generally, the spanning viewpoint of
\Cref{lem:first-integral-strict-gram} shows that this is structural.
Every square in \(I_2\) carries the common factor \((s+tv)^\top v\); at
\(s=e_1\) the generators are
\[
\widetilde q_{\alpha,t}(v)
=
(1+tv_1)^{\alpha_1}\prod_{j\ge2}(tv_j)^{\alpha_j}\,
\bigl(v_1+t\|v\|^2\bigr),
\qquad |\alpha|=d-2,
\]
whose only linear term is a multiple of \(v_1\).  For \(n\ge2\), the
coefficient vectors \(\widetilde a_\alpha(e_1,t)\) therefore never reach
the linear monomials \(v_2,\dots,v_n\): these coordinate directions lie in
the kernel of \(Q_2(e_1)\), and the spanning property fails.  Hence
\(Q_2(e_1)\succeq0\) but \(Q_2(e_1)\not\succ0\).  In fact, since each
\(t\mapsto\widetilde a_\alpha(e_1,t)\) is a polynomial of degree at most
\(d-1\),
\[
\operatorname{rank} Q_2(e_1)\;\le\; d\binom{n+d-3}{d-2},
\]
which is \(O(n^{d-2})\) for fixed \(d\), far below \(N=\Theta(n^d)\); by
orthogonal invariance, the same holds at every \(s\neq0\).  The contrast
with \(I_1\) is precisely the spanning mechanism: the squares in \(I_1\)
carry the prefactors \(v_i\) for \emph{every} coordinate \(i\), whereas the
squares in \(I_2\) share a single radial factor.  For \(n=1\), the
generators \(v_1(1+tv_1)^{d-1}\) do span \(v_1,\dots,v_1^d\) and
\(Q_2\succ0\), so the failure is a genuinely multivariate phenomenon.  The
strict Gram margin is therefore necessarily supplied by \(I_1\).

\begin{remark}[SoS-convexity and the spanning condition]
SoS-convexity of the regularization function is enough to show that its
Bregman remainder is SoS, but it does not by itself give a strict Gram margin.
Suppose a Bregman remainder admits the Gram representation
\[
R_s^\Phi(v)
=
\textbf{v}_d(v)^\top Q_\Phi(s)\textbf{v}_d(v),
\]
where
\[
Q_\Phi(s)
=
\int_0^1(1-t)
\sum_{\ell=1}^{L}
a_{\ell,s}(t)a_{\ell,s}(t)^\top\,dt
\succeq0.
\]
For this Gram matrix to be positive definite, the coefficient vectors
$a_{\ell,s}(t)$ must span the full coefficient space associated with
$\textbf{v}_d(v)$. Thus the additional requirement is that these coefficient vectors span the
full nonconstant monomial space; this is what upgrades an SoS representation
to a strict Gram representation. In the
Euclidean regularization considered here, the translated square terms
generated by the Hessian span all monomials of degrees $1,\ldots,d$,
which is precisely what yields the positive-definite Gram
representation in \Cref{lem:strict-gram}.
\end{remark}

\subsection{\texorpdfstring{Strict Gram Margins Beyond Euclidean Regularization}{Strict Gram Margins Beyond Euclidean Regularization}}
\label{sec: strict gram beyond}

The strict-Gram argument is not limited to the Euclidean norm.  What matters
is that the Bregman remainder can be written as an integral of polynomial
squares and that, together, those squares cover all monomial directions.  We
first record this principle and then apply it to several other regularizers,
including examples where the property holds and examples where it fails.

\begin{corollary}[Polynomial integral representations with a strict Gram margin]
\label{lem:polynomial-feature-gram}

Let \(z(v)\in\R^N\) be a fixed monomial vector.  Suppose that, for each
\(\ell\),
\[
    p_\ell(t,v)
    =
    a_\ell(t)^\top z(v),
\]
where \(a_\ell(t)\in\R^N\) is a vector-valued polynomial in \(t\),
\[
    a_\ell(t)
    =
    \sum_{k=0}^{L_\ell} b_{\ell,k}t^k,
    \qquad
    b_{\ell,k}\in\R^N.
\]
Suppose that the coefficient vectors satisfy the analogue of
\eqref{eq:annihilator-spanning-condition}, namely,
\begin{equation}
\label{eq:generalized-annihilator-spanning}
    \forall\,0\neq\xi\in\R^N,
    \qquad
    \exists\,\ell,\ 0\le k\le L_\ell
    \quad\text{s.t.}\quad
    \xi^\top b_{\ell,k}\neq0.
\end{equation}
Let \(c_\ell>0\), and let \(w:(0,1)\rightarrow\R_+\) be nonnegative and
strictly positive on a set of positive measure.  Then
\[
    Q
    :=
    \sum_\ell c_\ell
    \int_0^1
    w(t)\,
    a_\ell(t)a_\ell(t)^\top\,dt
\]
satisfies
$
    Q\succ0.
$
Equivalently, $
    \sum_\ell c_\ell
    \int_0^1
    w(t)\,p_\ell(t,v)^2\,dt
    =
    z(v)^\top Qz(v)$ 
has a strict Gram representation on \(z(v)\).

\end{corollary}

\begin{proof}
This follows from the same argument as
\Cref{lem:first-integral-strict-gram}.  Let \(0\neq\xi\in\R^N\).
By \eqref{eq:generalized-annihilator-spanning}, there exist
\(\ell\) and \(k\) such that
$
    \xi^\top b_{\ell,k}\neq0.
$
Hence
$
    \xi^\top a_\ell(t)
    =
    \sum_{k=0}^{L_\ell}
    \xi^\top b_{\ell,k}\,t^k
$
is not identically zero as a polynomial in \(t\).  Therefore,
$
    \int_0^1
    w(t)\bigl[\xi^\top a_\ell(t)\bigr]^2\,dt
    >0.
$
Since \(c_\ell>0\) and all the remaining terms are nonnegative,
$
    \xi^\top Q\xi
    =
    \sum_\ell c_\ell
    \int_0^1
    w(t)\bigl[\xi^\top a_\ell(t)\bigr]^2\,dt
    >0.
$
Thus \(Q\succ0\).
\end{proof}

\begin{remark}[Stability under SOS-convex additions]
Suppose that \(\Phi_0\) has a strict Bregman Gram representation on
\(\textbf{v}_d(v)\), and let \(\Psi\) be an SOS-convex polynomial of degree at
most \(2d\).  Then \(\Phi_0+\Psi\) also has a strict Bregman Gram
representation on \(\textbf{v}_d(v)\). 
Indeed, writing
$
R_{\Phi_0,s}(v)=\textbf{v}_d(v)^\top Q_0(s)\textbf{v}_d(v),
 Q_0(s)\succ0, $ $
R_{\Psi,s}(v)=\textbf{v}_d(v)^\top Q_\Psi(s)\textbf{v}_d(v),
  Q_\Psi(s)\succeq0,
$
we obtain
\[
R_{\Phi_0+\Psi,s}(v)
=
\textbf{v}_d(v)^\top\bigl(Q_0(s)+Q_\Psi(s)\bigr)\textbf{v}_d(v),
\]
where \(Q_0(s)+Q_\Psi(s)\succ0\).
\end{remark}

We now give three examples.  Weighted Euclidean regularization inherits the
Euclidean result by a change of coordinates.  Positive sums of linear forms can
also work when they generate enough polynomial directions.  By contrast,
coordinate-separable regularization misses all mixed directions and therefore
has no full strict Gram margin.

\paragraph{Weighted Euclidean radial regularizer.} 

Beyond the Euclidean regularizer, a natural
example is the weighted ellipsoidal regularization
$$
\Phi_W(x)
=
\frac{1}{2d}(x^\top Wx)^d,
 W\succ0.
$$
This is an invertible linear transformation of the Euclidean
regularizer. The induced transformations on the homogeneous monomial
spaces are also invertible, and therefore preserve the spanning
property and the existence of a positive-definite Gram representation.
The same conclusion holds for positive weighted sums of such
regularizers, for example
\[
\Phi(x)
=
\sum_{\ell=1}^{L}
\frac{\omega_\ell}{2d}
(x^\top M_\ell x)^d,
\qquad
\omega_\ell>0,
\]
provided that at least one $M_\ell$ is positive definite. More
generally, \Cref{lem:strict-gram} applies to any SoS-convex
regularizer whose Hessian-generated square terms satisfy the full monomial spanning condition (i.e., the Bregman remainder has a quantitative positive-definite Gram matrix).

\paragraph{Linear-form regularizers and Veronese spanning.} Consider
\[
\Phi_{\rm lf}(x)
=
\frac1{2d}
\sum_{\ell=1}^m
w_\ell(a_\ell^\top x)^{2d},
\qquad
w_\ell>0.
\]
The corresponding polynomials in the integral representation are
\[
p_\ell(t,v)
=
(a_\ell^\top v)(a_\ell^\top(s+tv))^{d-1}
=
\sum_{q=1}^{d}
\binom{d-1}{q-1}
(a_\ell^\top s)^{d-q}
t^{q-1}
(a_\ell^\top v)^q .
\]
Let \(\operatorname{coeff}((a_\ell^\top v)^q)\) denotes the coefficient\footnote{
Thus the
coefficient of \(t^{q-1}\) in \(p_\ell(t,v)\) is a scalar multiple of
$
\operatorname{coeff}((a_\ell^\top v)^q).
$
The scalar is
$
\binom{d-1}{q-1}(a_\ell^\top s)^{d-q}.
$
Therefore, provided \(a_\ell^\top s\neq0\), this scalar does not affect the
span.}
vector of the homogeneous degree-\(q\) polynomial \((a_\ell^\top v)^q\)
with respect to the monomial basis \(\{v^\beta:|\beta|=q\}\).  
Hence a sufficient condition for a strict Gram margin on the full nonconstant
monomial space is that, for every \(q=1,\ldots,d\),
\[
\operatorname{span}\{\operatorname{coeff}((a_\ell^\top v)^q):
\ell=1,\ldots,m\}
=
\mathcal H_{n,q}:=
\operatorname{span}\{e_{v^\beta}:|\beta|=q\}.
\]
Note that $\mathcal H_{n,q}$ denotes the coefficient space of homogeneous degree-\(q\) polynomials in \(v\).
This is exactly a Veronese spanning condition. In particular, for fixed s, a generic collection of sufficiently many linear forms can satisfy it.The Veronese embedding is precisely the map given by all degree-q monomials, turning degree-q polynomial data into linear coordinates.

This is exactly a Veronese spanning condition. 
In particular, for fixed $s$, a generic sufficient number of linear forms is
\[
m\ge \dim \mathcal H_{n,q}
=
\binom{n+d-1}{d}.
\]
In particular, for \(d=2\) it suffices generically to take
$
m\ge \binom{n+1}{2},
$
whereas for \(d=3\) it suffices generically to take
$
m\ge \binom{n+2}{3}.
$
The bound comes from the standard Veronese spanning condition 
\cite[Section~1.4]{Harris1992} \cite[Chapter~2]{Landsberg2012} and a proof is given in Appendix \ref{sec: Veronese} for completeness.

\paragraph{Failure of a strict margin for separable regularization.} Consider the separable regularizer
\[
\Phi_{\rm sep}(x)
=
\frac{1}{2d}\sum_{i=1}^{n}x_i^{2d}.
\]
It is SOS-convex, since
$
y^\top\nabla^2\Phi_{\rm sep}(x)y
=
(2d-1)\sum_{i=1}^{n}x_i^{2d-2}y_i^2
=
(2d-1)\sum_{i=1}^{n}
\left(x_i^{d-1}y_i\right)^2 .
$
Consequently, its Bregman remainder is SOS. However, the squared polynomials in the integral representation are only coordinate-separable. Indeed,
along the segment \(s+tv\), the Hessian representation involves terms of the
form
\[
    v_i(s_i+t v_i)^{d-1},
    \qquad i=1,\ldots,n.
\]
For each \(i\), this generates only pure powers of the single variable \(v_i\):
\[
v_i(s_i+t v_i)^{d-1}
=
\sum_{q=1}^{d}
\binom{d-1}{q-1}
s_i^{d-q}t^{q-1}v_i^q .
\]
Thus the coefficient directions generated by the separable regularizer are
contained in
\[
\operatorname{span}
\{e_{v_i^q}: i=1,\ldots,n,\ q=1,\ldots,d\},
\]
and no mixed monomial direction such as \(e_{v_iv_j}\), \(i\neq j\), is
generated. Hence the Bregman remainder cannot have a positive-definite Gram
representation on the full nonconstant monomial vector
$
\textbf{v}_d(v)=\bigl(v^\beta:1\le |\beta|\le d\bigr).
$
The failure is already visible on a coordinate line. Take \(s=e_1\) and
\(v=\tau e_2\). Then
\[
R_s^{\rm sep}(\tau e_2)
=
\frac{1}{2d}\tau^{2d}.
\]
If there existed a matrix \(Q(s)\succ0\) such that
$
R_s^{\rm sep}(v)
=
\textbf{v}_d(v)^\top Q(s)\textbf{v}_d(v),
$
then, since \(\textbf{v}_d(\tau e_2)\) contains the linear monomial \(\tau\), we would
have
\[
\frac{1}{2d}\tau^{2d}
=
R_s^{\rm sep}(\tau e_2)
\ge
\lambda_{\min}(Q(s))\,\|\textbf{v}_d(\tau e_2)\|^2
\ge
\lambda_{\min}(Q(s))\,\tau^2 .
\]
This is impossible as \(\tau\to0\) for \(d\ge2\). Therefore separable
regularization yields an SOS Bregman remainder, but it does not provide the
full strict Gram margin required by the perturbation argument.

    \begin{remark}[Coercivity does not imply a strict Gram margin]
Coercivity is not sufficient for a full strict Gram margin.  Indeed,
\[
\Phi_{\rm sep}(x)
=
\frac{1}{2d}\sum_{i=1}^n x_i^{2d}
=
\frac{1}{2d}\|x\|_{2d}^{2d}
\]
is coercive and SoS-convex, but it acts separately on each coordinate and
does not generate mixed monomial directions.  Hence, the key issue is not only
how strongly the regularizer grows, but whether it controls all polynomial
directions needed for a strict Gram representation.
\end{remark}

The conclusion of this section is that Euclidean regularization contributes a
strictly positive Gram core, while some other regularizers do not.  We now ask
whether this core is strong enough to dominate the rest of a higher-order
Taylor polynomial.

\section{\texorpdfstring{Sufficient Regularization, SoS Exactness, and Minimizer Recovery}{Sufficient Regularization, SoS Exactness, and Minimizer Recovery}}
\label{sec Main perturbation theorem}

This section turns the strict Gram margin into a global optimization result.
The idea is that a sufficiently large positive Gram core can dominate any
bounded lower-order terms that start at quadratic order.  We first quantify how the strict Gram margin absorbs bounded polynomial
perturbations, then apply this estimate to the same regularized polynomial
family written in symmetric-tensor Taylor form. Finally, we show that the
resulting strictness forces rank-one moment recovery of the minimizer.

\subsection{\texorpdfstring{Absorbing Bounded Polynomial Perturbations}{Absorbing Bounded Polynomial Perturbations}}

We quantify
how a strict Gram margin survives after adding a bounded polynomial
perturbation.

\begin{theorem}[SoS under bounded perturbations]
\label{thm:eventual-SoS}
{Let \(d\geq2\), fix \(s\in\R^n\setminus\{0\}\), and let
\(p\in\R[v]\) satisfy
\[
\deg p\leq 2d,\qquad
p(0)=0,\qquad
\nabla p(0)=0,
\qquad
\|\coeff(p)\|\le M
\]
for some \(M\ge0\).
{ Throughout this section, \(\textbf{v}_d(v)\) denotes the nonconstant monomial vector
$
\textbf{v}_d(v)=\bigl(v^\beta:1\leq |\beta|\leq d\bigr),
$
so the constant monomial is not included.}
For \(\sigma>0\), define
\[
q_\sigma(v)
:=
p(v)+\sigma R_{s,d}(v).
\]
There exists a constant \(C_{n,d}>0\), independent of \(s\), \(p\),
and \(\sigma\), such that whenever
\begin{equation}
\label{eq:fixed-threshold}
\sigma
\;>\;
\bar\sigma(s,M)
:=
\frac{C_{n,d}\,M}{c_{n,d}\,\min\{\|s\|^{2d-2},1\}},
\end{equation}
the polynomial \(q_\sigma\) admits a positive-definite Gram
representation with respect to \(\textbf{v}_d(v)\).  Consequently,
\(q_\sigma\in\SoS\) and
\(q_\sigma(v)\ge0\) with equality if and only if \(v=0\).}
\end{theorem}

\begin{proof}
By Lemma~\ref{lem:strict-gram},
\[
R_{s,d}(v)
=
\textbf{v}_d(v)^\top Q_R(s)\textbf{v}_d(v),
\qquad
Q_R(s)\succ0,
\]
with
$
\lmin(Q_R(s))
\geq
{c_{n,d}\min\{\|s\|^{2d-2},1\}.}
$
Because \(p\) has no constant or linear term and has degree at
most \(2d\), there exists a symmetric matrix \(Q_p\) such
that
\[
p(v)
=
\textbf{v}_d(v)^\top Q_p\textbf{v}_d(v).
\]
{The coefficient map from symmetric Gram matrices indexed by \(\textbf{v}_d(v)\) to polynomials
of degree at most \(2d\) with no constant or linear term is surjective.
Fixing a linear right inverse of this map,}
{whose operator norm defines the constant \(C_{n,d}\) in
\eqref{eq:fixed-threshold}, the coefficient bound gives}
\[
\|Q_p\|_2\leq C_{n,d}\,M.
\]
Therefore,
\[
q_\sigma(v)
=
\textbf{v}_d(v)^\top Q_\sigma \textbf{v}_d(v),
\qquad
Q_\sigma
:=
\sigma Q_R(s)+Q_p.
\]
By Weyl's inequality,
\begin{align*}
\lmin(Q_\sigma)
&\geq
\sigma\lmin(Q_R(s))
-
\|Q_p\|_2 \geq
{c_{n,d}\sigma\min\{\|s\|^{2d-2},1\}-C_{n,d}M,}
\end{align*}
{which is strictly positive precisely under the threshold condition
\eqref{eq:fixed-threshold}.} Hence
\(Q_\sigma\succ0\), and \(q_\sigma\) is SoS. Finally,
\[
q_\sigma(v)=0
\quad\Longleftrightarrow\quad
\textbf{v}_d(v)=0,
\]
because \(Q_\sigma\succ0\). The first homogeneous block of \(\textbf{v}_d(v)\)
is \([v]_1=v\), so \(\textbf{v}_d(v)=0\) if and only if \(v=0\). Thus
\(q_\sigma(v)>0\) for all \(v\neq0\).
\end{proof}

\begin{remark}[Why the perturbation starts at quadratic order]
\label{rem:quadratic-perturbation-space}
The assumptions
$p_\sigma(0)=0, 
\nabla p_\sigma(0)=0
$
are intrinsic to the present shifted strict-Gram argument.  Indeed, the
Bregman remainder satisfies
\[
R_{s,d}(0)=0,
\qquad
\nabla R_{s,d}(0)=0,
\]
and hence \(R_{s,d}(v)=O(\|v\|^2)\) as \(v\to0\).  Moreover, since the
Gram vector \(\textbf{v}_d(v)\) contains only nonconstant monomials, every
polynomial of the form
$
\textbf{v}_d(v)^\top Q\textbf{v}_d(v)
$
has zero constant and linear terms. 
Accordingly, the natural perturbation space is
\[
\mathcal P_{2d}^{(2)}
:=
\left\{
p\in\mathbb R[v]_{2d}:
p(0)=0,\ \nabla p(0)=0
\right\}.
\]
Every polynomial in \(\mathcal P_{2d}^{(2)}\) admits a symmetric Gram
representation on \(\textbf{v}_d(v)\).  The positive-definite Gram matrix of
\(R_{s,d}\) can therefore absorb any uniformly bounded signed Gram
perturbation in this space after sufficiently large scaling.
\end{remark}

\begin{corollary} 
\label{cor:eventual-SoS-family}
Let \(d\ge2\), and let \((p_\sigma)_{\sigma>0}\subset\R[v]\) and
\((s_\sigma)_{\sigma>0}\subset\R^n\setminus\{0\}\) satisfy
\[
\deg p_\sigma\le 2d,\qquad
p_\sigma(0)=0,\qquad
\nabla p_\sigma(0)=0,\qquad
\sup_\sigma\|\coeff(p_\sigma)\|<\infty,
\]
together with the scaling condition
\begin{equation}
\label{eq:scale-condition}
\sigma\|s_\sigma\|^{2d-2}\longrightarrow\infty .
\end{equation}
Then, for all sufficiently large \(\sigma\), the polynomial
\(q_\sigma=p_\sigma+\sigma R_{s_\sigma,d}\) admits a
positive-definite Gram representation with respect to
\(\textbf{v}_d(v)\); in particular \(q_\sigma\in\SoS\) and
\(q_\sigma(v)>0\) for \(v\neq0\).
\end{corollary}

\begin{proof}
Set \(M:=\sup_\sigma\|\coeff(p_\sigma)\|<\infty\).  By
\eqref{eq:scale-condition}, for all sufficiently large \(\sigma\) we
have \(\sigma\min\{\|s_\sigma\|^{2d-2},1\}>C_{n,d}M/c_{n,d}\),
so \Cref{thm:eventual-SoS} applies with the fixed pair
\((s,p)=(s_\sigma,p_\sigma)\) at each such \(\sigma\).
\end{proof}

\begin{remark}[Perturbed polynomial is nonconvex]
For every fixed $s$, the polynomial $R_{s,d}$ is SoS-convex. Indeed,
$
\nabla_v^2 R_{s,d}(v)
=
\nabla^2 \Phi_d(s+v),
$
and $\Phi_d(x)=\frac12\|x\|^{2d}$ is SoS-convex.
However, the perturbed polynomial
$
q_\sigma(v)
=
p_\sigma(v)+\sigma R_{s_\sigma,d}(v)
$
need not be convex, or SoS-convex. The proof of
Theorem~\ref{thm:eventual-SoS} uses the positive-definite Gram margin
of $R_{s_\sigma,d}$ in the polynomial SoS cone; it does not require
the Hessian of $q_\sigma$ to be positive semidefinite.
\end{remark}

\subsection{\texorpdfstring{Sufficiently Regularized Taylor Models}{Sufficiently Regularized Taylor Models}}

We now apply the perturbation result to the regularized polynomial
family of \eqref{eq:intro-regularized-model}, using the equivalent
symmetric-tensor Taylor form from \Cref{par:coefficient-matching}.
Its shifted coefficients remain bounded while the Gram margin grows.  The
assumption \(g\neq0\) is what gives the required scaling of the global
minimizer.

\begin{corollary}[SoS exactness of sufficiently regularized Taylor models]
\label{cor:SoS-taylor-model} 
Let \(d\geq2\), and consider
\[
m_\sigma(s)
=
\mathcal T_{2d-1}(s)
+
\frac{\sigma}{2d}\|s\|^{2d},
\]
where \(\mathcal T_{2d-1}\) is a fixed polynomial of degree at most
\(2d-1\). We assume that \(\|g\|\neq0\) and let \(s_\sigma^*\neq0\) be a global minimizer of
\(m_\sigma\). 
Then, for all sufficiently large \(\sigma\),
the shifted polynomial
$
q_\sigma(v)
:=
m_\sigma(s_\sigma^*+v)-m_\sigma(s_\sigma^*)
$
admits a positive-definite Gram representation with respect to
\(\textbf{v}_d(v)\) {(again excluding the constant monomial)} and
\[
q_\sigma\in\SoS,
\]
and
$
q_\sigma(v)=0$ if and only if $
v=0.
$
In particular, \(s_\sigma^*\) is the unique global minimizer of
\(m_\sigma\) for all sufficiently large \(\sigma\).

\end{corollary}

\begin{proof}
Since \(s_\sigma^*\) is an unconstrained global minimizer, it satisfies
$
\nabla\mathcal T_{2d-1}(s_\sigma^*)
+
\sigma\|s_\sigma^*\|^{2d-2}s_\sigma^*
=
0.
$
Expanding the shifted Taylor model around $s_\sigma^*$ and using stationarity, we obtain
$
q_\sigma(v)
=
p_\sigma(v)
+
{\frac{\sigma}{d}} R_{s_\sigma^*,d}(v),
$
{because \(R_{s,d}\) is defined for \(\Phi_d(x)=\frac12\|x\|^{2d}\), while
the Taylor model contains \(\frac{\sigma}{2d}\|s\|^{2d}\).}
where
\[
p_\sigma(v)
:=
\sum_{j=2}^{2d-1}
\frac{1}{j!}
\nabla^j\mathcal T_{2d-1}(s_\sigma^*)[v]^j.
\]
Clearly,
\[
\deg p_\sigma\leq2d-1,\qquad
p_\sigma(0)=0,\qquad
\nabla p_\sigma(0)=0.
\]
It remains to verify that the coefficients of \(p_\sigma\) are
uniformly bounded. For convenience, write \(T_2:=H\). For each
\(j=2,\ldots,2d-1\),
\[
\nabla^j\mathcal T_{2d-1}(s_\sigma^*)
=
\sum_{k=j}^{2d-1}
\frac{1}{(k-j)!}
T_k[(s_\sigma^*)]^{k-j}
\]
Since {\Cref{thm:general-bound-for-s} gives} \(s_\sigma^*\to0\), we have
\(\|s_\sigma^*\|\leq1\) for all sufficiently large \(\sigma\), and
therefore
\[
\left\|
\nabla^j\mathcal T_{2d-1}(s_\sigma^*)
\right\|
\leq
\sum_{k=j}^{2d-1}
\frac{\|T_k\|}{(k-j)!}.
\]
The right-hand side depends only on the fixed problem data and is
independent of \(\sigma\). Hence,
$
\sup_\sigma\|\operatorname{coeff}(p_\sigma)\|<\infty.
$ Using \Cref{thm:general-bound-for-s},  we have 
\[
\sigma\|s_\sigma^*\|^{2d-2}
=
\Theta\!\left(
\sigma^{1/(2d-1)}
\right)
\to\infty,
\]
which verifies the required scaling condition. {Equivalently,
\[
\frac{\sigma}{d}\|s_\sigma^*\|^{2d-2}\to\infty.
\]}
All the assumptions of \Cref{cor:eventual-SoS-family} are therefore
satisfied {with perturbation parameter \(\sigma/d\) in place of \(\sigma\)}, and the stated SoS and uniqueness conclusions follow.
\end{proof}

\begin{remark}[The degenerate case \(g=0\)]
The assumption \(g\neq0\) is essential for obtaining the universal
scaling
\[
\|s_\sigma^*\|
=
\Theta\!\left(\sigma^{-1/(2d-1)}\right).
\]
When \(g=0\), the scaling of a global minimizer is instead determined
by the first nonzero lower-order term of the Taylor polynomial, and
the global minimizer need not be unique.
For example, let \(n\geq2\), \(\mu>0\), and consider
\[
m_\sigma(s)
=
-\frac{\mu}{2}\|s\|^2
+
\frac{\sigma}{2d}\|s\|^{2d}.
\]
This corresponds to \(g=0\), \(H=-\mu I\), and
\(T_3=\cdots=T_{2d-1}=0\). 
Since \(r=0\) is a local maximum and
\(\varphi_\sigma(r)\to+\infty\) as \(r\to\infty\), the set of global
minimizers is
\[
\operatorname*{argmin}_{s\in\mathbb R^n}m_\sigma(s)
=
\left\{
s\in\mathbb R^n:
\|s\|
=
\left(\frac{\mu}{\sigma}\right)^{1/(2d-2)}
\right\}.
\]
Thus, there are infinitely many global minimizers, and every global
minimizer satisfies
\[
\|s_\sigma^*\|
=
\Theta\!\left(\sigma^{-1/(2d-2)}\right),
\]
rather than
\(\Theta\!\left(\sigma^{-1/(2d-1)}\right)\).
The
strict Gram-dominance condition \eqref{eq:scale-condition} of
\Cref{cor:eventual-SoS-family} is also not satisfied. For any chosen global
minimizer \(s_\sigma^*\), the shifted polynomial
$
q_\sigma(v)
=
m_\sigma(s_\sigma^*+v)-m_\sigma(s_\sigma^*)
$
has multiple zeros, corresponding to all other global minimizers, and
therefore cannot admit a positive-definite Gram representation with
respect to the full nonconstant monomial vector. Although \(q_\sigma\) may still be SoS and admit a positive-semidefinite Gram representation, it cannot admit a positive-definite Gram representation with respect to the full nonconstant monomial vector.
\end{remark}

\subsection{\texorpdfstring{Rank-One Moment Recovery and Exact Minimizer Extraction}{Rank-One Moment Recovery and Exact Minimizer Extraction}}

SoS exactness gives the correct global minimum value, but we also want the
minimizer itself.  The strict Gram representation provides this extra
information.  At an optimum, a positive-definite Gram matrix can pair to zero
with the positive-semidefinite moment matrix only if all shifted nonconstant
moments vanish.  The optimal moment matrix therefore collapses to rank one. 
In short,
\[
\boxed{
\text{strict Gram certificate}
\quad\Longrightarrow\quad
\text{unique rank-one moment solution}
\quad\Longrightarrow\quad
\text{exact minimizer recovery}.}
\]
The next theorem makes this statement precise.

\begin{theorem}[Rank-one moment recovery and exact minimizer extraction]
\label{thm:rank-one-moment-recovery}
Suppose that, for some \(d\geq 2\), a polynomial
\(m_\sigma\in\mathbb{R}[x]_{2d}\) has a unique global minimizer
\(s_\sigma^*\), and let
$
m_\sigma^*:=m_\sigma(s_\sigma^*).
$
Assume that the shifted polynomial
$
q_\sigma(v)
:=
m_\sigma(s_\sigma^*+v)-m_\sigma^*
$
admits the positive definite Gram representation
{with respect to the nonconstant monomial vector \(\textbf{v}_d(v)\)}
\begin{equation}
\label{eq:strict-shifted-gram}
q_\sigma(v)
=
\textbf{v}_d(v)^\top Q_\sigma \textbf{v}_d(v),
\qquad
Q_\sigma\succ0.
\end{equation} 
Then the order-\(d\) SoS representation
\begin{equation}
\label{eq:primal-SoS-relaxation}
\begin{aligned}
\rho_d^{\mathrm{SoS}}
:=
\sup_{\gamma\in\mathbb{R}}
\quad & \gamma\\
\text{\rm s.t.}\quad
& m_\sigma(x)-\gamma\in\Sigma[x]_{2d}
\end{aligned}
\end{equation}
is exact $
\rho_d^{\mathrm{SoS}}=m_\sigma^*.$ Moreover, the corresponding order-\(d\) moment relaxation
\begin{equation}
\label{eq:dual-moment-relaxation}
\begin{aligned}
\rho_d^{\mathrm{mom}}
:=
\inf_y
\quad & L_y(m_\sigma)\\
\text{\rm s.t.}\quad
& y_0=1,  M_d(y)\succeq0
\end{aligned}
\end{equation}
has a unique optimal truncated moment sequence \(y^*\), namely
$
y_\alpha^*=(s_\sigma^*)^\alpha, 
|\alpha|\leq 2d.
$
Consequently,
\[
M_d(y^*)
=
z_d(s_\sigma^*)
z_d(s_\sigma^*)^\top,
\]
where $
z_d(x)
=
\big(1,[x]_1^\top,\ldots,[x]_d^\top\big)^\top,
$
and hence
$
\operatorname{rank}M_d(y^*)
=
\operatorname{rank}M_{d-1}(y^*)
=
1.
$
In particular, the unique global minimizer is \textbf{recovered} exactly from
the first-order moments:
$
(s_\sigma^*)_i=y_{e_i}^*, 
i=1,\ldots,n.
$
\end{theorem}

\begin{proof}
\textbf{Primal feasible: Global minimum value attained.}
Since \(q_\sigma\) is SoS, the shifted polynomial
$
m_\sigma(x)-m_\sigma^*
=
q_\sigma(x-s_\sigma^*)
$
is also SoS by translations.
Therefore, \(\gamma=m_\sigma^*\) is feasible for
\eqref{eq:primal-SoS-relaxation}, and hence
$
\rho_d^{\mathrm{SoS}}\geq m_\sigma^*.
$
On the other hand, if
$
m_\sigma(x)-\gamma\in\Sigma[x]_{2d},
$
then \(m_\sigma(x)-\gamma\geq0\) for every \(x\), so necessarily
$
\gamma\leq\inf_x m_\sigma(x)=m_\sigma^*.
$
Thus
$
\rho_d^{\mathrm{SoS}}=m_\sigma^*.
$

\textbf{Dual feasible.} We next analyze the moment relaxation. The truncated moment sequence
of the Dirac measure at \(s_\sigma^*\),
$
y_\alpha=(s_\sigma^*)^\alpha,
 |\alpha|\leq2d,
$
is feasible for \eqref{eq:dual-moment-relaxation}, because
\[
M_d(y)
=
z_d(s_\sigma^*)
z_d(s_\sigma^*)^\top
\succeq0,
\qquad
y_0=1.
\]
Its objective value is
$
L_y(m_\sigma)=m_\sigma(s_\sigma^*)=m_\sigma^*.
$
Weak duality and the exactness of the SoS representation therefore yield
$
\rho_d^{\mathrm{mom}}
=
\rho_d^{\mathrm{SoS}}
=
m_\sigma^*.
$

\textbf{Global minimizer is recovered :} 
It remains to prove uniqueness and the rank-one property. 
{Let \(y\) be any optimal truncated moment sequence in the original
\(x\)-coordinates. Define the shifted truncated moment sequence
\(\widetilde y\) by
\[
L_{\widetilde y}(h(v))
=
L_y(h(x-s_\sigma^*)),
\qquad \deg h\leq 2d.
\]
Then,}
\[
\widetilde y_0=1,
\qquad
M_d(\widetilde y)\succeq0,
\qquad
L_{\widetilde y}(q_\sigma)=0.
\]
{Indeed, for every polynomial \(r\) of degree at most \(d\),
\[
L_{\widetilde y}(r(v)^2)
=
L_y(r(x-s_\sigma^*)^2)\geq0,
\]
so \(M_d(\widetilde y)\succeq0\); and optimality gives
\(L_{\widetilde y}(q_\sigma)=L_y(m_\sigma-m_\sigma^*)=0\).}
Partition the moment matrix according to the constant monomial and
the nonconstant monomial vector \(\textbf{v}_d(v)\):
\[
M_d(\widetilde y)
=
\begin{pmatrix}
1 & a^\top\\
a & M_d^+(\widetilde y)
\end{pmatrix},
\]
where \(M_d^+(\widetilde y)\) is the principal submatrix indexed by
all nonconstant monomials of degree at most \(d\). Since
\(M_d(\widetilde y)\succeq0\), we also have its  principal submatrix $
M_d^+(\widetilde y)\succeq0.
$
Using the strict Gram representation
\eqref{eq:strict-shifted-gram},
\[
\begin{aligned}
0
=
L_{\widetilde y}(q_\sigma)
=
L_{\widetilde y}
\bigl(\textbf{v}_d(v)^\top Q_\sigma \textbf{v}_d(v)\bigr)=
\left\langle
Q_\sigma,
M_d^+(\widetilde y)
\right\rangle = \trace (Q_\sigma
M_d^+(\widetilde y)).
\end{aligned}
\]
where $\left\langle
\cdot, \cdot
\right\rangle$ denotes the Frobenius inner product. {More details about the last two equalities are given in} \Cref{remark supplement 1}.
Because
\[
Q_\sigma\succ0
\qquad\text{and}\qquad
M_d^+(\widetilde y)\succeq0,
\]
the equality above implies
$
M_d^+(\widetilde y)=0.
$
 Indeed,
$
Q_\sigma^{1/2}
M_d^+(\widetilde y)
Q_\sigma^{1/2}
\succeq0
$
and its trace is zero, so it must be the zero matrix.
Hence
\[
M_d(\widetilde y)
=
\begin{pmatrix}
1&a^\top\\
a&0
\end{pmatrix}.
\]
The positive semidefiniteness of this matrix forces \(a=0\). For
example, every \(2\times2\) principal minor satisfies
$
\det
\begin{pmatrix}
1&a_i\\
a_i&0
\end{pmatrix}
=
-a_i^2\geq0,
$
and therefore \(a_i=0\) for every \(i\). Consequently,
$
M_d(\widetilde y)
=
\begin{pmatrix}
1&0\\
0&0
\end{pmatrix}.
$
{This also implies that all nonconstant moments vanish: the vector
\(a=0\) gives the moments of degrees \(1,\ldots,d\), while
\(M_d^+(\widetilde y)=0\) gives all moments of degrees \(2,\ldots,2d\), since
every monomial of degree between \(2\) and \(2d\) is a product of two
nonconstant monomials of degrees at most \(d\).}
Thus \(\widetilde y\) is \textbf{uniquely} the truncated moment sequence of
the Dirac measure at \(v=0\) (more details can be found in \Cref{remark supplement 2}.),
\[
\widetilde y_0=1,
\qquad
\widetilde y_\alpha=0
\quad
\text{for every }1\leq|\alpha|\leq2d.
\]
Translating back to the original variable \(x=v+s_\sigma^*\) gives
the unique optimal moment sequence
\[
y_\alpha^*=(s_\sigma^*)^\alpha,
\qquad |\alpha|\leq2d.
\]
Therefore,
$
M_d(y^*)
=
z_d(s_\sigma^*)
z_d(s_\sigma^*)^\top,
$
which has rank one. {The same conclusion holds for \(M_{d-1}(y^*)\), which is its principal
submatrix indexed by monomials of degree at most \(d-1\).}
\end{proof}

\begin{corollary}[Polynomial-size exact SDP and minimizer recovery]
\label{cor:polynomial-size-sdp}
Fix \(d\ge2\).  Under the assumptions of
\Cref{cor:SoS-taylor-model}, for all sufficiently large \(\sigma\), the
global minimum of \(m_\sigma\) is attained exactly by its degree-\(d\)
SoS/moment representation.  The corresponding moment matrix has order
$
    \binom{n+d}{d}.
$
Hence, for fixed \(d\), the size of the semidefinite formulation is
polynomial in \(n\).  Moreover, its unique optimal moment matrix has rank one
and its first-order moments recover the unique global minimizer.
\end{corollary}

\begin{proof}
Exactness follows from \Cref{cor:SoS-taylor-model}; uniqueness and rank-one
recovery follow from \Cref{thm:rank-one-moment-recovery}.  The matrix-order
formula is the number of monomials in \(n\) variables of degree at most \(d\).
\end{proof}

We have therefore completed the main chain of the paper:
\[
\text{large regularization}
\;\Longrightarrow\;
\text{strict Gram dominance}
\;\Longrightarrow\;
\text{exact SoS representation}
\;\Longrightarrow\;
\text{rank-one minimizer recovery}.
\]
The next section connects this unconstrained regularization picture with the
other classical globalization mechanism for Taylor models: trust regions.

\subsection{Further consequences of strict Gram representations}
\label{subsec:strict-gram-consequences}

A strict Gram representation has several implications beyond
SoS exactness and moment recovery. In particular, it guarantees
a unique global minimizer, nondegeneracy at the minimizer, and
a local Polyak--\L{}ojasiewicz (PL) inequality.
We first clarify the relationship between SoS exactness,
uniqueness, and strict Gram representations, including their
regularization thresholds. We then establish the local PL
property and show that an interior SoS-convex certificate
also implies a strict Gram representation.

\paragraph{SoS, unique global minimizers, and strict Gram}
A strict Gram representation implies both SoS exactness and
uniqueness of the global minimizer. Therefore,
\[
\sigma_{\mathrm{SoS}}\le \sigma_{\mathrm{strict}},
\qquad
\sigma_{\mathrm{uniq}}\le \sigma_{\mathrm{strict}},
\]
where $\sigma_{\mathrm{SoS}}$, $\sigma_{\mathrm{uniq}}$, and
$\sigma_{\mathrm{strict}}$ denote the smallest
regularization thresholds for SoS exactness, uniqueness of the
global minimizer, and the existence of a strict Gram
representation, respectively. Each denotes the infimum of parameters beyond which
the corresponding property holds for every larger regularization weight.
The three thresholds need not be attained at equality. 
Note that in general, neither SoS exactness nor uniqueness of the global minimizer implies the other. An SoS polynomial may have multiple
global minimizers, while a nonnegative polynomial with a unique
global minimizer need not be SoS. Consequently, the sufficient
regularization bound established above need not coincide with
the sharp threshold for either property.

\begin{example}
To illustrate this distinction, consider the quartically
regularized bivariate polynomial
\[
m_\sigma(x,y)
=x+\frac12x^2-\frac12y^2
+\frac{\sigma}{4}(x^2+y^2)^2,
\qquad \sigma>0.
\]
For every $\sigma>0$, completing squares gives
\[
m_\sigma(x,y)+\frac14+\frac{1}{4\sigma}
=
\left(x+\frac12\right)^2
+\frac{\sigma}{4}
\left(x^2+y^2-\frac1\sigma\right)^2.
\]
For $0<\sigma\le4$ the two squares vanish simultaneously, so
this identity gives the exact global minimum and an SoS certificate.
For $\sigma>4$, the shifted polynomial is also SoS by the strict
Gram representation stated below. For $0<\sigma<4$,
there are two global minimizers,
$
\left(-\frac12,\,
\pm\sqrt{\frac1\sigma-\frac14}\right).
$
At $\sigma=4$, they merge into the unique global minimizer
$(-1/2,0)$, where the Hessian is singular. 
For $\sigma>4$, the unique global minimizer is $(-u,0)$,
where $u\in(0,1/2)$ is the unique positive solution of
$\sigma u^3+u=1$. In this case, the shifted polynomial
admits a positive-definite Gram representation on the full
nonconstant monomial basis of degrees one and two. 
Therefore, the three sharp thresholds satisfy
\[
\sigma_{\mathrm{SoS}}=0
<
\sigma_{\mathrm{uniq}}
=\sigma_{\mathrm{strict}}=4.
\] 
\end{example}

\begin{remark}[Strict Gram and other stationary points]
A strict Gram representation guarantees a unique, nondegenerate
global minimizer, but does not exclude other stationary points
or nonglobal local minima. For example, consider
\[
f(x)=x^2\left((x-2)^2+\frac1{10}\right)
=x^4-4x^3+\frac{41}{10}x^2.
\]
It admits the strict Gram representation with $Q = \begin{pmatrix}
\frac{41}{10}&-2\\
-2&1
\end{pmatrix} \succ0$. 
Thus $x^*=0$ is the unique global minimizer.
Nevertheless, $f$ has two additional stationary points,
$
x_\pm=\frac{3\pm\sqrt{4/5}}{2},
$
where $x_-$ is a local maximum and $x_+$ is a
nonglobal local minimum.
Hence, strict Gram guarantees a unique
global minimizer, but does not imply that every stationary
point is globally optimal. 
The same phenomenon extends to saddle points in higher
dimensions: the polynomial
\[
F(x,y)=f(x)+y^2+x^2y^2+y^4
\]
admits a strict Gram representation on
$(x,y,x^2,xy,y^2)$, while $(x_-,0)$ is a saddle point
and $(x_+,0)$ is a nonglobal local minimizer.

\end{remark}

\paragraph{Local PL inequality from a strict Gram representation}

\begin{lemma}[Strict Gram implies local PL]
\label{lem:local-PL-gram}
Let $s^*$ be a global minimizer of a polynomial $m$, and suppose
that the shifted polynomial admits a strict Gram representation
\[
m(s^*+v)-m(s^*)
=\textbf{v}_d(v)^\top Q\textbf{v}_d(v),
\qquad Q\succ0.
\]
Then $s^*$ is a nondegenerate global minimizer. Moreover,
there exist a neighbourhood $U$ of $s^*$ and a constant $\mu>0$
such that
\[
m(s)-m(s^*)
\le \frac{1}{2\mu}\|\nabla m(s)\|^2,
\qquad s\in U.
\]
\end{lemma}

\begin{proof}
Partition the monomial vector and its Gram matrix as
\[
\textbf{v}_d(v)=
\begin{pmatrix}v\\ \textbf{v}_{2:d}(v)\end{pmatrix},
\qquad
Q=
\begin{pmatrix}
\frac12\nabla^2m(s^*) & Q_{12}\\
Q_{12}^\top & Q_{22}
\end{pmatrix},
\]
where $\textbf{v}_{2:d}$ contains the monomials of degrees
$2,\ldots,d$. The upper-left block follows by comparing the
quadratic terms in the Gram representation. Since $Q\succ0$,
we have $\nabla^2m(s^*)\succ0$. 
By continuity, there exists a sufficiently small convex
neighbourhood $U$ of $s^*$ on which
$\nabla^2m(s)\succeq\mu I$ for some $\mu>0$.
Thus $m$ is strongly convex on $U$, which gives the stated
local PL inequality.
\end{proof}

\begin{remark}[SoS, strict Gram, and local PL]
It is important to note that an SoS representation alone does not imply a local PL inequality, even when the global minimizer is unique. For example, $m(x)=x^4$ is SoS and has a unique global minimizer at $x^*=0$,
but it does not satisfy the standard local PL inequality near $x^*$.

In contrast, Lemma~\ref{lem:local-PL-gram} shows that a strict
Gram representation guarantees local PL. In fact, the key
condition in its proof is $\nabla^2m(s^*)\succ0$, which is
equivalent to positive definiteness of the degree-one principal
block of $Q$. Thus, full strict Gram is sufficient but not
necessary: a positive-semidefinite Gram representation with a
positive-definite degree-one block already implies local PL.
Conversely, local PL may hold even when $\nabla^2m(s^*)$
is singular.
\end{remark}

\paragraph{SoS-convexity and strict Gram representations.}

Ahmadi, Chaudhry and Zhang
\cite[Lemmas~2--3]{ahmadi2024higherorder}
show that $\|x\|^2+\|x\|^{2d}$ lies in the interior
of the SoS-convex cone, and use this property to construct
SoS-convex Taylor approximations.
This raises a natural question about the relationship
between interior SoS-convexity and the strict Gram
representations studied here. 
We first give an example showing that a strict Gram
representation may exist even when the polynomial is
nonconvex. We then show that an interior SoS-convex
certificate is sufficient for a strict Gram representation,
but not necessary. 

\begin{example}[Strict Gram before global convexification]
\label{ex:strict-before-convex}
Consider the bivariate quartic polynomial
\[
m_\sigma(x,y)
=\frac{137}{800}x+x^2+y^2+x^3
+\frac{\sigma}{4}(x^2+y^2)^2,
\qquad \sigma>0.
\]
It has a nonzero linear term and is locally strongly
convex at the origin, since
$\nabla^2m_\sigma(0,0)=2I$. 
For every $\sigma>0$, the shifted polynomial
$m_\sigma-m_\sigma^*$ is SoS, by the classical
result on nonnegative bivariate quartics.
Moreover, at $\sigma=5/4$, the unique global
minimizer is $s^*=(-1/10,0)$, and the shifted
polynomial admits a positive-definite Gram
representation on all nonconstant monomials
of degrees one and two.
Nevertheless, $m_{5/4}$ is not globally convex.
Indeed, direct calculation gives the exact
convexity threshold
\[
\sigma_{\mathrm{convex}}=\frac32, \qquad \text{and} \qquad \sigma_{\mathrm{SoSconvex}}=\frac32.
\]
as convexity and SoS-convexity coincide for
bivariate quartics \cite{ahmadi2013complete}.
Thus, this example has
\[
\sigma_{\mathrm{SoS}}=0,\qquad 
\sigma_{\mathrm{strict}}\le\frac54
<
\sigma_{\mathrm{convex}}
=\sigma_{\mathrm{SoSconvex}}=\frac32.
\]
In particular, strict Gram exactness can occur
before global convexification, even when the
polynomial is locally strongly convex at the
regularization centre.
\end{example}

 Helton and Nie
\cite[Lemma~8]{helton2010semidefinite}
showed that a polynomial with an SoS Hessian and a
stationary zero is SoS, while Ahmadi and Parrilo
\cite[Theorem~3.1]{ahmadi2013complete}
established the corresponding first-order SoS
characterization of SoS-convexity.
The following lemma records a stronger conclusion
when the Hessian SoS representation has a
positive-definite Gram matrix.

\begin{lemma}[Interior SoS-convexity implies strict Gram]
\label{lem:interior-sosconvex-strict}
Let $h$ be a strictly SoS-convex polynomial of degree at most $2d$, namely
\[
y^\top\nabla^2h(x)y
=\mathbf w(x,y)^\top M\mathbf w(x,y),
\qquad M\succ0,
\]
where
$
\mathbf w(x,y)
=
\bigl(y_i x^\alpha\bigr)_{
1\le i\le n,\ |\alpha|\le d-1}.
$
Then $h$ has a unique global minimizer $s^*$, and the
shifted polynomial $h(s^*+v)-h(s^*)$ admits a strict
Gram representation on $\textbf{v}_d(v)$.
\end{lemma}

\begin{proof}
By stationarity and Taylor's integral formula,
\[
h(s^*+v)-h(s^*)
=
\int_0^1(1-t)\,
v^\top\nabla^2h(s^*+tv)v\,dt.
\]
Substituting the assumed Hessian representation,
the integrand is a positive-definite quadratic
form in the polynomials
\[
v_i(s^*+tv)^\alpha,
\qquad 1\le i\le n,\quad |\alpha|\le d-1.
\]
Let $Q$ denote the resulting Gram matrix on
$\textbf{v}_d(v)$. If $Q$ were singular, there
would exist a nonzero vector $\xi$ orthogonal
to the coefficient vectors of all these
polynomials for almost every $t\in(0,1)$.
Since these coefficient vectors are polynomial
in $t$, the same orthogonality would hold
at $t=1$. 
However, the polynomials
$v_i(s^*+v)^\alpha$ span all polynomials
with zero constant term and degree at most $d$:
translation preserves the polynomial space of
degree at most $d-1$, and multiplication by
$v_i$ generates every nonconstant monomial.
Thus $\xi=0$, a contradiction.
Consequently, $Q\succ0$.
\end{proof}

\section{\texorpdfstring{Trust Regions and Regularization}{Trust Regions and Regularization}}
\label{sec:TRplus-formulation}

The previous section used regularization.  We now compare it with a
trust-region constraint.  The main observation is that, when the trust-region
solution lies on the boundary, its multiplier selects a matching regularization
parameter.  With this choice the two shifted polynomial gaps are the same.
This lets the SoS certificate pass from one formulation to the other.  We then
show that small trust-region radii correspond to large regularization, and we
treat the interior cubic case separately.

We have therefore completed the main chain of the paper:
\[
\text{large regularization}
\;\Longrightarrow\;
\text{strict Gram dominance}
\;\Longrightarrow\;
\text{exact SoS representation}
\;\Longrightarrow\;
\text{rank-one minimizer recovery}.
\]
The next section connects this unconstrained regularization picture with the
other classical globalization mechanism for Taylor models: trust regions.

\subsection{\texorpdfstring{Matching Trust-Region and Regularized Models}{Matching Trust-Region and Regularized Models}}
\label{subsec:why-TRplus}

We first make the boundary correspondence exact.  The \(\mathrm{TR}^+\)
formulation keeps the same feasible ball but uses a higher even radial term in
the certificate.  The benefit is that, at an active solution, this certificate
can be matched exactly with the regularized Taylor model.
Recall that
\begin{equation}
    \mathcal T_{2d-1}(s)=f_0+g^\top s+\sum_{j=2}^{2d-1}\frac1{j!}T_j[s]^j,
    \label{Td}
\end{equation}
with \(d\ge2\). The standard trust-region problem is
\begin{equation}
s_{\Delta}^* \in \argmin_{\|s\|\le\Delta}\ \mathcal T_{2d-1}(s),
\tag{TR}
\label{eq:TR}
\end{equation}
with standard Lagrangian
\begin{equation}
\label{Lagrangian TR}
   L_{\rm TR}(s,\lambda)
   =
{\mathcal T_{2d-1}(s)}
    +\frac{\lambda}{2}\bigl(\|s\|^2-\Delta^2\bigr),
   \qquad
   \lambda\ge0.
  \tag{Lagrangian-TR}
\end{equation}
Its active-boundary stationarity condition is
$\nabla \mathcal T_{2d-1}(s_\Delta^*)+\lambda_\Delta s_\Delta^*=0,
   \|s_\Delta^*\|=\Delta .$

We define a new trust-region formulation, denoted by \(\mathrm{TR}^+\), which
uses the same feasible ball but replaces the quadratic localizer by a higher
radial one.  Equivalently, one may raise the radial constraint to degree
\(2d\):
\begin{equation}
{
\|s\|^2\le \Delta^2
\qquad\Longleftrightarrow\qquad
\|s\|^{2d}\le \Delta^{2d}.
}
\tag{\(\mathrm{TR}^+\)}
\label{eq:TR plus}
\end{equation}
For \(d\ge2\), we define the \(\mathrm{TR}^+\) Lagrangian polynomial {(a normalized higher-radial Lagrangian, i.e.\ a global certificate polynomial)} by
\begin{equation}
\mathcal{L}_{\Delta,\lambda}(s)
:=
\mathcal T_{2d-1}(s)
+
\frac{\lambda}{2d\,\Delta^{2d-2}}\bigl(\|s\|^{2d}-\Delta^{2d}\bigr).
\tag{Lagrangian-\(\mathrm{TR}^+\)}
\label{eq:TRplus}
\end{equation}
Clearly, the trust-region constraint \(\|s\|^2\le\Delta^2\) and
the higher radial constraint \(\|s\|^{2d}\le\Delta^{2d}\) define the same ball.
Moreover, under the normalization
$
\frac{1}{2d\,\Delta^{2d-2}}
\bigl(\|s\|^{2d}-\Delta^{2d}\bigr),
$
they have the same active-boundary stationarity equation
$
   \nabla \mathcal T_{2d-1}(s^*)+\lambda s^*=0.
$
Indeed, the normalization is chosen so that, on \(\|s\|=\Delta\),
\[
\nabla
\left[
\frac{1}{2d\,\Delta^{2d-2}}
\bigl(\|s\|^{2d}-\Delta^{2d}\bigr)
\right]_{\|s\|=\Delta}
= s .
\]
We say that the \(\mathrm{TR}^+\) certificate \eqref{eq:TRplus} is
\emph{exact} at an active trust-region minimizer \(s_\Delta^*\), with multiplier
\(\lambda_\Delta\), if
\begin{eqnarray}
    \mathcal{L}_{\Delta,\lambda_\Delta}(s)
-
\mathcal{L}_{\Delta,\lambda_\Delta}(s_\Delta^*)
\in\Sigma[s]_{2d}.
\label{lag gap}
\end{eqnarray}
where $\lambda_\Delta$ represents the Lagrangian multiplier corresponding to $s_\Delta^*$.

Thus, when the constraint is active, \(\mathrm{TR}\) and \(\mathrm{TR}^+\) have
the same KKT equation, but they play different roles.  The standard
\(\mathrm{TR}\) formulation is the constrained optimization problem, whereas
\(\mathrm{TR}^+\) uses the same feasible ball to provide a global SOS
certificate.  The third related object is the unconstrained regularized model
\begin{equation}
s_\sigma^* \in \argmin_{s \in \mathbb R^n}
m_{\sigma}(s),
\qquad
m_{\sigma}(s):=
\mathcal T_{2d-1}(s)+\frac{\sigma}{2d}\|s\|^{2d}.
\tag{Regularized model}
\label{eq:regularized-model}
\end{equation}
The next proposition makes precise the shared active-boundary KKT geometry of
\(\mathrm{TR}\) and \(\mathrm{TR}^+\), and the matching relation with the
regularized model when the radius is chosen consistently.

\begin{proposition}[Matched trust-region and regularized formulations]
\label{prop:TR-TRplus-regularized-matching}
Let \(d\ge2\) and \(\Delta>0\), and let \(s_\Delta^*\) be a
global minimizer of the trust-region problem \eqref{eq:TR} at which
the constraint is active, \(\|s_\Delta^*\|=\Delta\).  Then the
following hold.
\begin{enumerate}
\item \(s_\Delta^*\) satisfies the active-boundary KKT equation
$
   \nabla \mathcal T_{2d-1}(s_\Delta^*)+\lambda_\Delta s_\Delta^*=0
$
for some multiplier \(\lambda_\Delta\ge0\).
\item \textbf{(Matching.)}  Set
\(\sigma:=\lambda_\Delta/\Delta^{2d-2}\).  Then the matched
\(\mathrm{TR}^+\) Lagrangian gap coincides, as a polynomial in
\(s\), with the shifted regularized gap:
\[
\mathcal{L}_{\Delta,\lambda_\Delta}(s)
-
\mathcal{L}_{\Delta,\lambda_\Delta}(s_\Delta^*)
=
m_{\sigma}(s)-m_{\sigma}(s_\Delta^*),
\]
with \(m_\sigma\) as in \eqref{eq:regularized-model}.  Consequently,
the \(\mathrm{TR}^+\) Lagrangian gap is SOS exact at
\((\Delta,\lambda_\Delta)\) if and only if the shifted regularized
polynomial \(m_{\sigma}(s)-m_{\sigma}(s_\Delta^*)\) is SOS exact.
\item \textbf{(Converse matching.)}  If \(s_\sigma^*\neq0\) is a
global minimizer of \(m_\sigma\) in \eqref{eq:regularized-model},
then, with \(\Delta:=\|s_\sigma^*\|\), the point \(s_\sigma^*\) is
a global minimizer of the trust-region problem
\(\min_{\|s\|\le\Delta}\mathcal T_{2d-1}(s)\), the constraint is active,
and items~1--2 hold with multiplier
\(\lambda_\Delta=\sigma\Delta^{2d-2}>0\).
\end{enumerate}
\end{proposition}

\begin{proof}
(1) At an active minimizer, the gradient of the constraint
\(\|s\|^2-\Delta^2\) is \(2s_\Delta^*\neq0\), so the linear
independence constraint qualification holds on the boundary and the
KKT conditions produce a multiplier \(\lambda_\Delta\ge0\) with
\(\nabla \mathcal T_{2d-1}(s_\Delta^*)+\lambda_\Delta s_\Delta^*=0\).

(2) Since \(\|s_\Delta^*\|=\Delta\), the localizer vanishes at
\(s_\Delta^*\), and with \(\sigma=\lambda_\Delta/\Delta^{2d-2}\),
\[
\frac{\lambda_\Delta}{2d\,\Delta^{2d-2}}
\bigl(\|s\|^{2d}-\Delta^{2d}\bigr)
=
\frac{\sigma}{2d}
\bigl(\|s\|^{2d}-\|s_\Delta^*\|^{2d}\bigr).
\]
Adding \(\mathcal T_{2d-1}(s)-\mathcal T_{2d-1}(s_\Delta^*)\) to both sides gives the
polynomial identity of item~2; the SOS equivalence is immediate from
the identity.

(3) Since \(s_\sigma^*\) is a global minimizer of \(m_\sigma\), for
every \(s\) with \(\|s\|\le\Delta=\|s_\sigma^*\|\),
\[
\mathcal T_{2d-1}(s)-\mathcal T_{2d-1}(s_\sigma^*)
\ge
\frac{\sigma}{2d}
\bigl(\Delta^{2d}-\|s\|^{2d}\bigr)
\ge0,
\]
so \(s_\sigma^*\) is a global minimizer of the trust-region problem
with radius \(\Delta\), attained on the boundary because
\(s_\sigma^*\neq0\).  The first-order condition
\(\nabla \mathcal T_{2d-1}(s_\sigma^*)
+\sigma\|s_\sigma^*\|^{2d-2}s_\sigma^*=0\) is the KKT equation of
item~1 with \(\lambda_\Delta=\sigma\Delta^{2d-2}\), and item~2 applies
verbatim.
\end{proof}

\subsection{\texorpdfstring{Small Trust Regions and Large Regularization}{Small Trust Regions and Large Regularization}}
\label{subsec:small-radius-main}

The matching formula has a simple consequence: shrinking the trust-region
radius increases the matching regularization parameter.  Thus the
small-radius regime is the trust-region version of the large-regularization
regime studied in \Cref{sec Main perturbation theorem}.

\Cref{prop:TR-TRplus-regularized-matching} identifies shrinking the
trust-region radius with increasing the regularization weight: the
matched parameter is \(\sigma=\lambda_\Delta/\Delta^{2d-2}\), and, by the
boundary expansion of \Cref{lem:small-radius-boundary} below, the
active multiplier satisfies
\(\lambda_\Delta\Delta\to\|g\|\) as \(\Delta\downarrow0\), so
that
\[
    \sigma
    =
    \frac{\lambda_\Delta}{\Delta^{2d-2}}
    =
    \frac{\|g\|}{\Delta^{2d-1}}\bigl(1+o(1)\bigr)
    \longrightarrow\infty .
\]
Consequently, when a trust-region scheme keeps shrinking the radius at
a nonstationary point, it is implicitly driving the matched
regularized model into the large-\(\sigma\) regime of
\Cref{sec Main perturbation theorem}.  The next lemma gives the first
half of this picture: an explicit radius threshold, depending only on
the coefficients of the higher-order Taylor terms, below which the
global trust-region minimizer is forced onto the boundary.

For the small-radius estimate below, let
\[
    \Lambda_j
    :=
    \max_{\|v_1\|=\cdots=\|v_j\|=1}
    |T_j[v_1,\ldots,v_j]|
\]
denote the operator norm of the symmetric tensor \(T_j\)
(and hence \(\Lambda_2=\|H\|_2\)).

The next lemma shows that, when the linear term is nonzero, sufficiently small
trust-region subproblems are forced to attain their global minimizer on the
boundary.  The radius threshold is controlled explicitly by the tensor norms of
the higher-order terms.

\begin{lemma}[Global minimizer on the boundary for small radius]
\label{lem:small-radius-boundary}
Let \(\mathcal T_{2d-1}\) be defined as in \eqref{Td}, and assume \(g\neq0\).  
Define
\[
   A_T
   :=
   \sum_{j=2}^{2d-1}
   \frac{\Lambda_j}{(j-1)!}.
\]
where the constant \(A_T\) depends only on the coefficients \(H,T_3,\dotsc,T_{2d-1}\) of \(\mathcal T_{2d-1}\).
{%
Let
\begin{eqnarray}
\label{radiu small}
   \bar\Delta_{\rm bd}
   :=
   \begin{cases}
      1, & A_T=0,\\[2pt]
      \min\!\left\{1,\dfrac{\|g\|}{2A_T}\right\}, & A_T>0.
   \end{cases}
\end{eqnarray}
Then, for every \(0<\Delta\le \bar\Delta_{\rm bd}\), every global minimizer
$
   s_\Delta^*\in\argmin_{\|s\|\le\Delta}\mathcal T_{2d-1}(s)
$
lies on the boundary,
$
   \|s_\Delta^*\|=\Delta.
$
After possibly decreasing \(\bar\Delta_{\rm bd}\), the minimizer is unique and satisfies
$
   s_\Delta^*
   =
   -\Delta\dfrac{g}{\|g\|}
   +
   O(\Delta^2).
$}
\end{lemma}

\begin{proof}
Write $\mathcal T_{2d-1}(s)=f_0+g^\top s+E(s)$, where
$\nabla E(s)=O(\|s\|)$ on the unit ball.  By the definition of $A_T$,
$\|\nabla E(s)\|\le A_T\|s\|$ for $\|s\|\le1$.
An interior minimizer with $\|s\|<\Delta\le\bar\Delta_{\rm bd}$
would satisfy $g+\nabla E(s)=0$, and hence
$\|g\|\le A_T\Delta\le\|g\|/2$, a contradiction.

For uniqueness and the expansion, set $s=\Delta y$ and consider
$F_\Delta(y)=[\mathcal T_{2d-1}(\Delta y)-f_0]/\Delta$ on
$\|y\|\le1$.  Then $F_\Delta(y)=g^\top y+O(\Delta)$ uniformly
in $C^2$.  The limiting linear problem has the unique solution
$y_0=-g/\|g\|$ on the sphere, so all minimizers approach $y_0$.
In a sufficiently small spherical neighborhood of $y_0$, the restriction
of $g^\top y$ has a positive-definite tangential Hessian.  This remains
true for $F_\Delta$ when $\Delta$ is small, giving a unique minimizer
$y_\Delta=y_0+O(\Delta)$ and thus the claimed expansion.
\end{proof}

\begin{proposition}[Small-radius SoS exactness]\label{thm:small-radius-TRplus}
Assume $g\neq0$.  For all sufficiently small $\Delta>0$, the global
minimizer $s_\Delta^*$ of the trust-region problem is unique and lies
on the boundary.  Its matched $\mathrm{TR}^+$ Lagrangian gap admits
a strict Gram representation in $v=s-s_\Delta^*$.  Consequently,
the degree-$d$ moment relaxation is exact and recovers $s_\Delta^*$
from its unique rank-one optimal moment matrix.
\end{proposition}
\begin{proof}
By \Cref{lem:small-radius-boundary},
$s_\Delta^*=-\Delta g/\|g\|+O(\Delta^2)$.
The active KKT equation and $\nabla\mathcal T_{2d-1}(s_\Delta^*)=g+O(\Delta)$
therefore give $\lambda_\Delta\Delta\to\|g\|$, so
$\lambda_\Delta\to\infty$ as $\Delta\downarrow0$.
By \Cref{prop:TR-TRplus-regularized-matching}, the shifted $\mathrm{TR}^+$
gap equals the shifted model with
$\sigma=\lambda_\Delta/\Delta^{2d-2}$.
Its Taylor perturbation has no constant or linear term and has uniformly
bounded coefficients as $s_\Delta^*\to0$.  The Gram margin supplied by
$\frac{\sigma}{d}R_{s_\Delta^*,d}$ is at least
$\frac{c_{n,d}}d\sigma\Delta^{2d-2}
=\frac{c_{n,d}}d\lambda_\Delta\to\infty$.
Thus \Cref{thm:eventual-SoS} gives a strict Gram matrix for sufficiently
small $\Delta$, and \Cref{thm:rank-one-moment-recovery} gives rank-one
recovery.
\end{proof}

\begin{remark}[Use within adaptive regularization and trust-region frameworks]
\label{rem:ARp-use}
The polynomial subproblems considered here are closely related to those
used in adaptive regularization and trust-region methods
\cite[Chs.~3--4]{cartis2022evaluation}; see also
\cite{birgin2017worst,cartis2020concise}. The standard AR$p$ method
uses a $p+1$ power, whereas our SoS analysis requires an even
polynomial degree; the trust-region certificate likewise uses the
matched higher-radial $\mathrm{TR}^+$ formulation.
At a nonstationary iterate, sufficiently large $\sigma_k$, or a
sufficiently small trust-region radius $\Delta_k$, yields an exact
SoS certificate and rank-one recovery of the model minimizer.
Under the usual smoothness and acceptance-ratio assumptions, a
sufficiently large regularization parameter also ensures successful
AR$p$ iterations. Thus a certified inner solve can be incorporated
into an adaptive scheme without changing its acceptance mechanism.
A complete convergence and complexity analysis, including uniform
thresholds, inexact SDP solutions and the cost of inner solves,
is beyond the scope of this paper.
\end{remark}

\subsection{\texorpdfstring{Interior Solutions in the Cubic Case}{Interior Solutions in the Cubic Case}}
\label{subsec:interior-local}

The matching argument above concerns boundary solutions.  If a trust-region
minimizer lies strictly inside the ball, the constraint is locally inactive,
so the point is simply an unconstrained local minimizer of the polynomial.
For cubic models, this case can be handled by the semidefinite local-minimum
results of Ahmadi and Zhang.

\begin{proposition}[Interior trust-region local minima]
\label{prop:interior-local-min}
Let \(\mathcal T_3\) be a cubic polynomial and let \(\Delta>0\).  If \(s^*\) is a local
minimizer of
$
   \min_{\|s\|\le\Delta}\mathcal T_3(s)
$
and satisfies \(\|s^*\|<\Delta\), then \(s^*\) is an unconstrained local
minimizer of \(\mathcal T_3\).  Conversely, any unconstrained local minimizer of \(\mathcal T_3\)
that lies strictly inside the ball is a local minimizer of the trust-region
problem.
\end{proposition}

\begin{proof}
If \(\|s^*\|<\Delta\), then there exists a neighbourhood \(U\) of \(s^*\) such
that
$
   U\subset\{s:\|s\|<\Delta\}.
$
Hence constrained local minimality on the trust-region ball is equivalent to
unconstrained local minimality on \(U\).  The converse follows from the same
argument in reverse.
\end{proof}

\begin{remark}[Ahmadi--Zhang local-minimum oracle]
\label{rem:ahmadi-zhang-oracle}
For cubic polynomials, the interior case can be handled by the semidefinite
machinery of Ahmadi and Zhang~\cite{ahmadi2022complexity}.  They show that
local minimality of a point of a cubic polynomial admits an efficiently
checkable necessary and sufficient condition, and that a local minimizer can be
found by solving semidefinite programs whose matrix size is linear in the
number of variables.  Thus, in the cubic case \(d=2\), the interior branch can
be treated by a cubic local-minimum SDP oracle.  The boundary/small-radius
branch is handled instead by the \(\mathrm{TR}^+\) certificate in
\Cref{thm:small-radius-TRplus}.
\end{remark}

The two trust-region cases therefore fit together cleanly.  On a sufficiently
small active boundary, the multiplier matches the problem to the
large-regularization regime and gives the same strict Gram/SoS certificate.  In the
cubic interior case, one instead uses an unconstrained local-minimum oracle.
This completes the connection between the strict-Gram theory above and
trust-region globalization.

\section{\texorpdfstring{Conclusion}{Conclusion}}

\label{sec:conclusion}

This paper studies when the polynomial subproblems arising in higher-order
optimization become exactly solvable through a single SoS/moment semidefinite
representation.  Our main result shows that, under a nonzero linear term and
sufficient Euclidean regularization, the nonnegative polynomial obtained by
shifting a higher-order Taylor model by its global minimum is SoS.  For every
fixed degree, the resulting SDP has polynomial size in the dimension, and its
rank-one moment solution recovers the unique global minimizer.  Importantly,
this conclusion does not require the regularized polynomial to be convex or
SoS-convex.

The proof identifies the structural reason behind this phenomenon.  The
Bregman remainder of the Euclidean radial regularizer has a strict Gram
representation at every nonzero point, with a strict Gram margin.  This
margin eventually dominates the bounded lower-order Taylor terms.  The same
viewpoint also explains a genuine boundary: some regularizers inherit
the strict Gram property, whereas coordinate-separable regularization does
not control the mixed polynomial directions required by the argument.  Thus
not only the amount and degree of regularization, but also its algebraic
structure, determine SoS exactness.

The trust-region results provide a second interpretation.  At an active
boundary solution, the multiplier defines a matched radial regularization,
and the small-radius regime corresponds precisely to the large-regularization
regime of the main theorem.  For cubic models, the complementary interior
case is also tractable through the semidefinite characterization of cubic
local minima.  

Several broader questions remain open.  The first is how far the tractable
regime identified in this paper can be extended beyond sufficient
regularization.  Our results show that sufficiently strong regularization can
turn a structured nonconvex higher-order polynomial subproblem into an exact
SoS/SDP problem of polynomial size for every fixed degree.  A natural next
question is whether the same phenomenon persists under substantially weaker
regularization.  General polynomial optimization is of course
intractable, so such a result cannot hold without additional structure; the
challenge is to identify how much structure is sufficient for exact global
optimization to remain possible.

A second fundamental issue is scalability.  Although the SoS/moment
formulation obtained here has polynomial size for fixed degree, standard SDP
representations still grow rapidly with the dimension and can become
expensive in practice.  There are, however, reasons for optimism.  The proof
developed in this paper does more than establish SoS exactness: it exposes
specific algebraic structure in the Gram representation generated by the
regularization term.  Related low-order models already show that exploiting
such structure can lead to substantially smaller tight SDP formulations,
including linear-size formulations in certain cubic--quartic settings
\cite{zhou2025tight,zhu2023cubic}.  This suggests that symmetry, sparsity,
low-rank tensor structure, and the particular Gram structure revealed here
may be used to construct much smaller exact or approximate formulations.
Turning these structural observations into scalable global solvers, and
ultimately incorporating them into globally convergent higher-order methods,
is an important direction for future research.

More broadly, the results suggest that the computational difficulty of
higher-order polynomial models is governed not only by their degree, but by
the algebraic structure created by regularization.  Understanding and
exploiting this structure may provide a route from worst-case polynomial
hardness to globally solvable and practically scalable subclasses.

\paragraph{Declaration of generative AI and AI-assisted technologies}
During the preparation of this work, the authors used large language
model systems (Claude, Anthropic; Codex models, OpenAI; Qwen, Alibaba) to assist in
reviewing mathematical proofs for potential errors or gaps. The authors
independently assessed the feedback and take full responsibility for
the correctness and content of the publication.

\paragraph{Funding.}
 The research of Coralia Cartis
and Wenqi Zhu was also supported by the Engineering and Physical Sciences
Research Council under grant EP/Y028872/1, \emph{Mathematical Foundations
of Intelligence: An ``Erlangen Programme'' for AI}.

\medskip
\noindent\textbf{Competing interests.}
The authors have no relevant financial or non-financial interests to disclose.

\medskip
\noindent\textbf{Data availability.}
No datasets were generated or analysed during the current study.

\appendix

\section{Full Statements of the Cited SoS-Convexity Results}
\label{app:cited-sos-results}

For completeness, and to make the citations in
\Cref{rem:standard-reference} precise, we reproduce the two results of
Ahmadi--Parrilo and Helton--Nie in the form in which they are used in
\Cref{sec The SoS property for the Bregman}.

\begin{theorem}[{Ahmadi--Parrilo \cite[Thm.~3.1]{ahmadi2013complete}}]
\label{thm:app-ahmadi-parrilo}
Let \(f\in\R[x]\) be a polynomial.  The following are equivalent:
\begin{enumerate}[label=\textup{(\roman*)}]
\item the Hessian \(\nabla^2f(x)\) is an SoS matrix, i.e.\
      \(y^\top\nabla^2f(x)\,y\) is a sum of squares in \((x,y)\)
      \textup{(}\(f\) is SoS-convex\textup{)};
\item the first-order divergence
      \(f(y)-f(x)-\nabla f(x)^\top(y-x)\) is a sum of squares in
      \((x,y)\);
\item the midpoint gap
      \(\tfrac12 f(x)+\tfrac12 f(y)-f\bigl(\tfrac{x+y}{2}\bigr)\)
      is a sum of squares in \((x,y)\).
\end{enumerate}
In particular, SoS-convexity of \(f\) implies that every
first-order Bregman divergence of \(f\) is SoS, jointly in the base point
and the increment.
\end{theorem}

\begin{lemma}[{Helton--Nie \cite[Lem.~8]{helton2010semidefinite}}]
\label{lem:app-helton-nie}
Let \(f\in\R[x]\) and suppose there exists a point \(u\in\R^n\) such
that
\[
    f(u)=0,
    \qquad
    \nabla f(u)=0,
\]
and the Hessian \(\nabla^2 f(x)\) is an SoS matrix.  Then \(f\) is a
sum of squares.
\end{lemma}

The point \(u\) in \Cref{lem:app-helton-nie} is the second-order
contact point of \(f\): \(f\) vanishes and is stationary at \(u\), so
that in the shifted variable \(v=x-u\) the polynomial has no constant
or linear part.  In our setting the lemma is applied to
\(f(v)=R_{s,d}(v)\) with contact point \(u=0\): indeed
\(R_{s,d}(0)=0\), \(\nabla_vR_{s,d}(0)=0\), and
\(\nabla_v^2R_{s,d}(v)=\nabla^2\Phi_d(s+v)\) is an SoS matrix by the
computation in the proof of \Cref{lem:bregman-SoS}.  Either result
therefore yields the qualitative SoS property of the Bregman remainder
\(R_{s,d}\).  Neither result is quantitative: the strict Gram margin
\eqref{eq:margin}, on which \Cref{thm:eventual-SoS} rests, is
established independently in \Cref{lem:strict-gram}.

\section{Generic Veronese spanning}
\label{sec: Veronese}

This is the usual Veronese embedding/spanning viewpoint; see, e.g.,
\cite[Section~1.4]{Harris1992} for the Veronese embedding in algebraic
geometry, or \cite[Chapter~2]{Landsberg2012} for the symmetric tensor
interpretation.
\begin{lemma}[Generic Veronese spanning]
Let
$
\mathcal H_{n,q}
:=
\operatorname{span}\{e_{v^\beta}:|\beta|=q\}
$
denote the coefficient space of homogeneous degree-\(q\) polynomials in
\(v\in\mathbb R^n\).  For \(a\in\mathbb R^n\), define
$
\nu_q(a)
:=
\operatorname{coeff}\bigl((a^\top v)^q\bigr)
\in \mathcal H_{n,q}.
$
Let \(d\ge1\), and suppose
\[
m\ge
\max_{1\le q\le d}\dim \mathcal H_{n,q}
=
\max_{1\le q\le d}\binom{n+q-1}{q}
=
\binom{n+d-1}{d}.
\]
Then, for a generic collection of vectors
$
a_1,\ldots,a_m\in\mathbb R^n,
$
we have
\[
\operatorname{span}\{\nu_q(a_\ell):\ell=1,\ldots,m\}
=
\mathcal H_{n,q},
\qquad
q=1,\ldots,d.
\]
Equivalently, the Veronese matrices (i.e., multivariate homogeneous Vandermonde matrix)
\[
M_q
=
\bigl[\nu_q(a_1)\ \cdots\ \nu_q(a_m)\bigr]
\in
\mathbb R^{\binom{n+q-1}{q}\times m}
\]
have full row rank for every \(q=1,\ldots,d\).  Moreover, for fixed
\(s\neq0\), the same generic choice can be made so that
\(a_\ell^\top s\neq0\) for every \(\ell\).
\end{lemma}

\begin{proof}
First note that
\[
\dim \mathcal H_{n,q}
=
\binom{n+q-1}{q}, \qquad \max_{1\le q\le d}\dim\mathcal H_{n,q}
=
\binom{n+d-1}{d}.
\]
as this quantity is nondecreasing in \(q\). 
Fix \(q\in\{1,\ldots,d\}\).  The entries of \(M_q\) are polynomial functions
of the vectors \(a_1,\ldots,a_m\).  Therefore the condition
$
\operatorname{rank} M_q < \dim\mathcal H_{n,q}
$
is an algebraic condition, given by the vanishing of all maximal minors of
\(M_q\).  It remains to show that these minors are not all identically zero. 
Let
$
N_q:=\dim \mathcal H_{n,q}.
$
The homogeneous monomials \(\{a^\beta:|\beta|=q\}\) are linearly independent
polynomial functions of \(a\).  Hence one can choose points
\[
b_1,\ldots,b_{N_q}\in\mathbb R^n
\]
such that the homogeneous Vandermonde matrix
$
\bigl[b_j^\beta\bigr]_{|\beta|=q,\ j=1,\ldots,N_q}
$
is nonsingular.  Since \(\nu_q(b_j)\) differs from
\((b_j^\beta)_{|\beta|=q}\) only by nonzero multinomial factors on the rows,
the matrix
\[
\bigl[\nu_q(b_1)\ \cdots\ \nu_q(b_{N_q})\bigr]
\]
also has full rank \(N_q\).  Thus at least one maximal minor of \(M_q\) is
not the zero polynomial.
Consequently, for each fixed \(q\), the set of vectors
\((a_1,\ldots,a_m)\) for which \(M_q\) has full row rank is the complement
of a proper algebraic set.  Taking the finite intersection over
\(q=1,\ldots,d\), we obtain a generic set on which all matrices
\(M_q\) have full row rank simultaneously.
Finally, for fixed \(s\neq0\), the conditions
$
a_\ell^\top s\neq0,
 \ell=1,\ldots,m,
$
exclude only finitely many proper hyperplanes.  Intersecting with their
complement preserves genericity.  This proves the claim.
\end{proof}

 \section{Scaling of Global Minimizers under Large Regularization}

\begin{theorem}[Scaling of the global minimizer]
\label{thm:general-bound-for-s}
Assume that \(\|g\|\neq0\), and let \(s_\sigma^*\) satisfy
\[
\nabla m_\sigma(s_\sigma^*)=0,
\qquad
m_\sigma(s_\sigma^*)\leq m_\sigma(0).
\]
Then there exists a constant \(\bar{\sigma}>0\), depending only on the
fixed problem data, such that, for every \(\sigma\geq\bar{\sigma}\),
\begin{equation}
\label{eq:general-bound-for-s}
\left(\frac{\|g\|}{2\sigma}\right)^{\frac{1}{2d-1}}
\leq
\|s_\sigma^*\|
\leq
\left(\frac{3d\|g\|}{\sigma}\right)^{\frac{1}{2d-1}}.
\end{equation}
Consequently,
\[
\|s_\sigma^*\|
=
\Theta\!\left(\sigma^{-\frac{1}{2d-1}}\right)
\qquad\text{as }\sigma\to\infty.
\]
\end{theorem}

\begin{proof}
Let
\[
r_\sigma:=\|s_\sigma^*\|,
\qquad
a_j:=\frac{\Lambda_j}{j!},
\qquad
b_j:=\frac{\Lambda_j}{(j-1)!},
\]
where \(\Lambda_j\) is the tensor operator norm introduced before \Cref{lem:small-radius-boundary}.
Since \(m_\sigma(s_\sigma^*)\leq m_\sigma(0)\), we have
\[
0
\geq
m_\sigma(s_\sigma^*)-m_\sigma(0)
\geq
-\|g\|r_\sigma
-\sum_{j=2}^{2d-1}a_jr_\sigma^j
+\frac{\sigma}{2d}r_\sigma^{2d}.
\]
We first show that \(r_\sigma\to0\). Let
$
C:=\|g\|+\sum_{j=2}^{2d-1}a_j.
$
If \(r_\sigma\geq1\), then
$
\frac{\sigma}{2d}r_\sigma^{2d}
\leq
Cr_\sigma^{2d-1},
$
and hence
$
r_\sigma\leq\frac{2dC}{\sigma},
$
which contradicts \(r_\sigma\geq1\) for all sufficiently large
\(\sigma\). Therefore \(r_\sigma<1\) eventually, and
$
\frac{\sigma}{2d}r_\sigma^{2d}
\leq
Cr_\sigma,
$
so
$
r_\sigma^{2d-1}
\leq
\frac{2dC}{\sigma}.
$
Thus \(r_\sigma\to0\). 
Since \(r_\sigma\to0\), for all sufficiently large \(\sigma\),
$
\sum_{j=2}^{2d-1}a_jr_\sigma^j
\leq
\frac{\|g\|}{2}r_\sigma.
$
Using again \(m_\sigma(s_\sigma^*)\leq m_\sigma(0)\), we obtain
\[
\frac{\sigma}{2d}r_\sigma^{2d}
\leq
\frac{3\|g\|}{2}r_\sigma.
\]
Since \(\|g\|\neq0\), stationarity implies \(r_\sigma>0\), and hence
\[
r_\sigma^{2d-1}
\leq
\frac{3d\|g\|}{\sigma}.
\]
This proves the upper bound in \eqref{eq:general-bound-for-s}.

For the lower bound, the stationarity condition gives
\[
g
+
\sum_{j=2}^{2d-1}
\frac{1}{(j-1)!}
T_j[(s_\sigma^*)]^{j-1}
+
\sigma r_\sigma^{2d-2}s_\sigma^*
=
0.
\]
Taking norms and using the reverse triangle inequality yields
\[
\sigma r_\sigma^{2d-1}
\geq
\|g\|
-
\sum_{j=2}^{2d-1}
b_jr_\sigma^{j-1}.
\]
Since \(r_\sigma\to0\), for all sufficiently large \(\sigma\),
$
\sum_{j=2}^{2d-1}
b_jr_\sigma^{j-1}
\leq
\frac{\|g\|}{2}.
$
Therefore,
$
\sigma r_\sigma^{2d-1}
\geq
\frac{\|g\|}{2},
$
which gives
$
r_\sigma
\geq
\left(\frac{\|g\|}{2\sigma}\right)^{\frac{1}{2d-1}}.
$
Combining the two bounds proves
$
\|s_\sigma^*\|
=
\Theta\!\left(\sigma^{-\frac{1}{2d-1}}\right).
$
\end{proof}




\section{More Discussions on Moment Sequence}

\begin{remark}
Let
$
\mathcal A_d
:=
\left\{
\alpha\in\mathbb N^n:
1\le |\alpha|\le d
\right\}, 
|\alpha|:=\alpha_1+\cdots+\alpha_n.
$
Then
\[
\textbf{v}_d(v)
=
\bigl(v^\alpha\bigr)_{\alpha\in\mathcal A_d},
\qquad
v^\alpha:=v_1^{\alpha_1}\cdots v_n^{\alpha_n},
\]
is the vector of all nonconstant monomials in \(v\) of degree at most
\(d\). Then, 
\[
q_\sigma(v)
=
\sum_{\alpha,\beta\in\mathcal A_d}
(Q_\sigma)_{\alpha\beta}v^{\alpha+\beta},
\]
while
$
\bigl(M_d^+(\widetilde y)\bigr)_{\alpha\beta}
=
\widetilde y_{\alpha+\beta}
=
L_{\widetilde y}(v^{\alpha+\beta}).
$
Hence, 
$
L_{\widetilde y}(q_\sigma)
=
\sum_{\alpha,\beta}
(Q_\sigma)_{\alpha\beta}
\widetilde y_{\alpha+\beta}
=
\left\langle
Q_\sigma,M_d^+(\widetilde y)
\right\rangle.
$
\label{remark supplement 1}
\end{remark}

\begin{remark}
The unique optimal truncated moment
sequence is the moment sequence of the Dirac measure at \(v=0\). More
precisely,
\[
\widetilde y_0=1,
\qquad
\widetilde y_\alpha=0
\quad
\text{for every }1\leq |\alpha|\leq 2d.
\]
Indeed, for the Dirac measure \(\delta_0\),
\[
\int_{\mathbb R^n} v^\alpha\,d\delta_0(v)
=
\begin{cases}
1, & \alpha=0,\\
0, & |\alpha|\geq1.
\end{cases}
\]
Hence \(\widetilde y\) coincides with the truncated moment sequence of
\(\delta_0\). The uniqueness follows from the preceding argument,
which shows that any optimal shifted moment matrix must equal
$
M_d(\widetilde y)
=
\begin{pmatrix}
1&0\\
0&0
\end{pmatrix}.
$
For
$
x=v+s_\sigma^*,
$
under this translation, the Dirac measure at \(v=0\) is mapped to the
Dirac measure at \(x=s_\sigma^*\). Therefore, the unique optimal
truncated moment sequence in the original coordinates is
\[
y_\alpha^*
=
\int_{\mathbb R^n}x^\alpha\,d\delta_{s_\sigma^*}(x)
=
(s_\sigma^*)^\alpha,
\qquad
|\alpha|\leq2d.
\]
Consequently, for every multi-index
\(\alpha,\beta\in\mathbb N^n\) satisfying
\(|\alpha|,|\beta|\leq d\), the corresponding entry of the moment
matrix is
\[
\bigl(M_d(y^*)\bigr)_{\alpha,\beta}
=
y_{\alpha+\beta}^*
=
(s_\sigma^*)^{\alpha+\beta}
=
(s_\sigma^*)^\alpha(s_\sigma^*)^\beta.
\]
Hence
\[
M_d(y^*)
=
z_d(s_\sigma^*)
z_d(s_\sigma^*)^\top,
\]
where
$
z_d(x)
=
\bigl(x^\alpha\bigr)_{|\alpha|\leq d}
=
\bigl(1,[x]_1^\top,\ldots,[x]_d^\top\bigr)^\top.
$
Since \(z_d(s_\sigma^*)\neq0\), we have 
$
\operatorname{rank}M_d(y^*)=1.
$
Likewise,
\[
M_{d-1}(y^*)
=
z_{\leq d-1}(s_\sigma^*)
z_{\leq d-1}(s_\sigma^*)^\top,
\]
which also have
$
\operatorname{rank}M_{d-1}(y^*)=1.
$
Thus the flat rank condition holds:
\[
\operatorname{rank}M_d(y^*)
=
\operatorname{rank}M_{d-1}(y^*)
=
1.
\]
Finally, for the standard unit multi-index \(e_i\),
$
y_{e_i}^*
=
(s_\sigma^*)^{e_i}
=
(s_\sigma^*)_i,
 i=1,\ldots,n.
$
Therefore, the unique global minimizer can be recovered directly from
the first-order moments:
$
s_\sigma^*
=
\bigl(
y_{e_1}^*,\ldots,y_{e_n}^*
\bigr)^\top.
$
Finally, the first-order moments satisfy
$
y_{e_i}^*=(s_\sigma^*)_i,
$
which gives the claimed exact extraction formula.
    \label{remark supplement 2}
\end{remark}

\scriptsize{
\bibliographystyle{plain}
\bibliography{sample}
}

\end{document}